\documentclass[11pt]{article}
\usepackage{amstext,amssymb,amsmath,amsbsy}
\usepackage{hyperref}
\usepackage{amscd}
\usepackage{amsfonts}
\usepackage{indentfirst}
\usepackage{verbatim}
\usepackage{amsmath}
\usepackage{amsthm}
\usepackage{enumerate}
\usepackage{graphicx}
\usepackage{color, soul}
\usepackage[OT1]{fontenc}
\usepackage[latin1]{inputenc}
\usepackage[english]{babel}
\usepackage{amssymb}
\usepackage{subfig}
\usepackage{algorithm}
\usepackage{algorithmic}
\usepackage{mathtools}

\newtheorem{remark}{Remark}[section]

\newtheorem*{Assumption*}{Assumption}

\newtheorem*{problem*}{Problem}
\numberwithin{equation}{section}
\input{tcilatex}
\begin{document}

\title{Carleman contraction mapping for a 1D inverse scattering problem with
experimental time-dependent data}
\author{Thuy T. Le\footnotemark[1] \and Micheal V. Klibanov\thanks{%
Department of Mathematics and Statistics, University of North Carolina at
Charlotte, Charlotte, NC, 28223, USA, \texttt{tle55@uncc.edu,
mklibanv@uncc.edu (corresponding author) loc.nguyen@uncc.edu}.} \and Loc H.
Nguyen\footnotemark[1] \and Anders Sullivan\thanks{%
US Army Research Laboratory, 2800 Powder Mill Road, Adelphi, MD 20783-1197,
USA, (anders.j.sullivan.civ@mail.mil, lam.h.nguyen2civ@mail.mil)} \and Lam
Nguyen\footnotemark[2] }
\date{}
\maketitle

\begin{abstract}
It is shown that the contraction mapping principle with the involvement of a
Carleman Weight Function works for a Coefficient Inverse Problem for a 1D
hyperbolic equation. Using a Carleman estimate, the global convergence of
the corresponding numerical method is established. Numerical studies for
both computationally simulated and experimentally collected data are
presented. The experimental part is concerned with the problem of computing
dielectric constants of explosive-like targets in the standoff mode using
severely underdetermined data.
\end{abstract}

\noindent {\textit{Key words}: one-dimensional wave equation, Carleman
estimate, iterative method, contraction principle, numerical solution,
experimental data, global convergence. }

\noindent{\textit{AMS subject classification}: 35R25, 35R30}

\section{Introduction}

\label{sec:1}

We consider in this paper a 1D Coefficient Inverse Problem (CIP) for a
hyperbolic PDE and construct a globally convergent numerical method for it.
Unlike all previous works of this group on the convexification method for
CIPs, which are cited below in this section, we construct in this paper a
sequence of linearized boundary value problems with overdetermined boundary
data. Each of these problems is solved using a weighted version of the
Quasi-Reversibility Method (QRM). The weight is the Carleman Weight Function
(CWF), i.e. this is the function, which is involved as the weight in the
Carleman estimate for the corresponding PDE operator, see, e.g. \cite%
{Carl,BK,BY,KT,Ksurvey,KL,LRS,Yam} for Carleman estimates. Thus, we call
this method \textquotedblleft Carleman Quasi-Reversibility Method" (CQRM).
The key convergence estimate for this sequence is similar with the one of
the conventional contraction mapping principle, see the first item of
Remarks 8.1 in section 8. This explains the title of our paper. That
estimate implies the global convergence of that sequence to the true
solution of our CIP. Furthermore, that estimate implies a rapid convergence.
As a result, computations are done here in \emph{real time}, also, see
Remark 11.1 in section 11 as well as section 12. The real time computation
is obviously important for our target application, which is described below.
The numerical method of this paper can be considered as the \emph{second
generation of the convexification} \emph{method} for CIPs.

The convexification method ends up with the minimization of a weighted
globally strictly convex Tikhonov functional, in which the weight is the
CWF. The convexification has been developed by our research group since the
first two publications in 1995 \cite{KlibIous} and 1997 \cite{Klib97}. Only
the theory \ was presented in those two originating papers. The start for
numerical studies of the convexification was given in the paper \cite{BAK}.
Thus, in follow up publications \cite%
{KlibKol1,KlibKol2,KLNSN,Khoa1,Khoa2,Khoa3,KlibKol1,KlibKol2,KLZ1,KLZ2,KLZ3,KlibSAR1,KlibSAR2,KLNSN,LN2,SKN1,SKN2}
the theoritical results for the convexification are combined with the
numerical ones. In particular, some of these papers work with experimentally
collected backscattering data \cite{Khoa2,Khoa3,KlibKol3,KlibSAR1,KlibSAR2}.
We also refer to the recently published book of Klibanov and Li \cite{KL},
in which many of these results are described.

The CIP of our paper has a direct application in the problem of the standoff
detection and identification of antipersonnel land mines and improvised
explosive devices (IEDs). Thus, in the numerical part of this publication,
we present results of the numerical performance of our method for both
computationally simulated and experimentally collected data for targets
mimicking antipersonnel land mines and IEDs. The experimental data were
collected in the field, which is more challenging than the case of the data
collected in a laboratory \cite{Khoa2,Khoa3,KlibKol3,KlibSAR2}. The data
used in this paper were collected by the forward looking radar of the US
Army Research Laboratory \cite{NguyenWong}. Since these data were described
in the previous publications of our research group \cite%
{Karch,KlibLoc,KlibKol1,KlibKol2,KLNSN,Kuzh1,Kuzh2,SKN2}, then we do not
describe them here. Those previous publications were treating these
experimental data by several different numerical methods for 1D CIPs. We
also note that the CIP of our paper was applied in \cite{KlibSAR1,KlibSAR2}
to the nonlinear synthetic aperture (SAR) imaging, including the cases of
experimentally collected data.

Our goal is to compute approximate values of dielectric constants of
targets. We point out that our experimental data are \emph{severely
underdetermined}. Indeed, any target is a 3D object. On the other hand, we
have only one experimentally measured time resolved curve for each target.
Therefore, we can compute only a sort of an average value of the dielectric
constant of each target. This explains the reason of our mathematical
modeling here by a 1D hyperbolic PDE rather than by its 3D analog. This
mathematical model, along with its frequency domain analog, was considered
in all above cited previous publications of our research group, which were
working with the experimental data of this paper by a variety of globally
convergent inversion techniques.

An estimate of the dielectric constant of a target cannot differentiate with
a good certainty between an explosive and a non explosive. However, we hope
that these estimates, combined with current targets classification
procedures, might potentially result in a decrease of the false alarm rate
for those targets. We believe therefore that our results for experimental
data might help to address an important application to the standoff
detection and identification of antipersonnel land mines and IEDs.

Given a numerical method for a CIP, we call this method globally convergent
if:

\begin{enumerate}
\item There is a theorem claiming that this method delivers at least one
point in a sufficiently small neighborhood of the true solution of that CIP\
without relying on a good initial guess about this solution.

\item This theorem is confirmed numerically.
\end{enumerate}

CIPs are both highly nonlinear and ill-posed. As a result, conventional
least squares cost functionals for them are non convex and typically have
many local minima and ravines, see, e.g. \cite{Scales} for a numerical
example of this phenomenon. Conventional numerical methods for CIPs are
based on minimizations of those functionals, see, e.g. \cite%
{Chavent,Gonch1,Gonch2}. The goal of the convexification is to avoid local
minima and ravines and to provide accurate and reliable solutions of CIPs.

In the convexification, a CWF is a part of a certain least squares cost
functional $J_{\lambda }$, where $\lambda >1$ is the parameter of the CWF.
The presence of the CWF ensures that, with a proper choice of the parameter $%
\lambda ,$ the functional $J_{\lambda }$ is strictly convex on a certain
bounded convex set $S$ of an arbitrary diameter $d>0$ in a Hilbert space. It
was proven later in \cite{BAK} that the existence and uniqueness of the
minimizer of $J_{\lambda }$ on $S$ are also guaranteed. Furthermore, the
gradient projection method of the minimization of $J_{\lambda }$ on $S$
converges globally to that minimizer if starting at an arbitrary point of $S$%
. In addition, if the level of noise $\delta >0$ in the data tends to zero,
then the gradient projection method delivers a sequence, which converges to
the true solution. Thus, since $d>0$ is an arbitrary number, then this is
the global convergence, as long as numerical studies confirm this theory.
Furthermore, it was recently established in \cite{KlibSAR2,LN2} that the
gradient descent method of the minimization of $J_{\lambda }$ on $S$ also
converges globally to the true solution.

First, we come up with the same boundary value problem (BVP) for a nonlinear
and nonlocal hyperbolic PDE as we did in the convexification method for the
same 1D CIP in \cite{KLNSN}. The solution of this BVP immediately provides
the target unknown coefficient. However, when solving this BVP, rather than
constructing a globally strictly convex cost functional for this BVP, as it
was done in \cite{KLNSN} and other works on the convexification, we
construct a sequence of linearized BVPs. Just as the original BVP, each BVP\
of this sequence involves a hyperbolic operator with non local terms and
overdetermined boundary conditions.

Because of non local terms and the overdetermined boundary conditions, we
solve each BVP of this sequence by the new version of the QRM, the Carleman
Quasi-Reversibility Method (CQRM). Indeed, the classical QRM solves
ill-posed/overdetermined BVPs for linear PDEs, see \cite{LL} for the
originating work as well as, e.g. \cite{KlibAPNUM,KL} for updates.
Convergence of the QRM-solution to the true solution is proven via an
appropriate Carleman estimates \cite{KlibAPNUM,KL}. However, a CWF is not
used in the optimizing quadratic functional of the conventional QRM of \cite%
{KlibAPNUM,KL}. Unlike this, a CWF is involved in CQRM. It is exactly the
presence of the CWF, which enables us to derive the above mentioned
convergence estimate, which is an analog of the one of the contraction
principle. This allows us, in turn to establish convergence rate of the
sequence of CQRM-solutions to the true solution of that BVP. The latter
result ensures that our method is a globally convergent one.

Various versions of the second generation of the convexification method
involving CQRM were previously published in \cite{Baud1,Baud2,LK,
NguyenKlibanov:preprint2021,LN1}. In these publications, globally convergent
numerical methods for some nonlinear inverse problems are presented. In \cite%
{Baud1,Baud2} CIPs for hyperbolic PDEs in $\mathbb{R}^{n}$ are considered.
The main difference between \cite{Baud1,Baud2} and our work is that it is
assumed in \cite{Baud1,Baud2} that one of initial conditions is not
vanishing everywhere in the closed domain of interest. In other words,
papers \cite{Baud1,Baud2} work exactly in the framework of the
Bukhgeim-Klibanov method, see \cite{BukhKlib} for the originating work on
this method and, e.g. \cite{BK,BY,Klib92,KT,Ksurvey,KL,Yam} for some follow
up publications. On the other hand, in all above cited works on the
convexification, including this one, either the initial condition is the $%
\delta -$function, or a similar requirement holds for the Helmholtz
equation. The only exception is the paper \cite{KLZ3}, which also works
within the framework of the method of \cite{BukhKlib}. We refer to papers 
\cite{Hoop,Korpela} for some numerical methods for CIPs with the
Dirichlet-to-Neumann map data. In these works the number $m$ of free
variables in the data exceeds the number $n$ of free variables in the
unknown coefficient, $m>n.$ In our paper $m=n=1.$ Also, $m=n$ in all other
above cited works on the convexification.

There is a classical Gelfand-Levitan method \cite{Levitan} for solutions of
1D CIPs for some hyperbolic PDEs. This method does not rely on the
optimization, and, therefore, avoids the phenomenon of local minima. It
reduces the original CIP to a linear integral equation of the second kind.
This equation is called \textquotedblleft Gelfand-Levitan equation" (GL).
However, the questions about uniqueness and stability of the solution of GL\
for the case of noisy data is open, see, e.g. Lemma 2.4 in the book \cite[%
Chapter 2]{Rom}. This lemma is valid only in the case of noiseless data.
However, the realistic data are always noisy. In addition, it was
demonstrated numerically in \cite{Karch} that GL cannot work well for
exactly the same experimental data as the ones we use in the current paper.
On the other hand, it was shown in \cite{KlibKol1,KlibKol2,SKN2} that the
convexification works well with these data. The same is true for the method
of the current paper.

Uniqueness and Lipschitz stability theorems of the CIP considered here are
well known. Indeed, it was shown in, e.g. \cite{SKN2} that, using the same
change of variables as the one considered in section 2 below, one can reduce
our CIP to a similar CIP for the equation $v_{tt}=v_{yy}+r\left( y\right)
v,y\in \mathbb{R}$ with the unknown coefficient $r\left( y\right) .$ We
refer to \cite[Theorem 2.6 of Chapter 2]{Rom} for the Lipschitz stability
estimate for the latter CIP. In addition, both uniqueness and Lipschitz
stability results for our CIP actually follow from Theorem 8.1 below as well
as from the convergence analysis in two other works of our research group on
the convexification method for this CIP \cite{KLNSN,SKN2}.

This paper is arranged as follows. In section 2 we state both forward and
inverse problems. In section 3 we derive a nonlinear BVP with non local
terms. In section 4 we describe iterations of our method to solve that BVP.
In section 5 we formulate the Carleman estimate for the principal part of
the PDE operator of that BVP. In section 6 we prove the strong convexity of
a functional of section 5 on an appropriate bounded set. In section 7 we
formulate two methods for finding the unique minimizer of that functional:
gradient descent method and gradient projection method and prove the global
convergence to the minimizer for each of them. In section 8 we establish the
contraction mapping property and prove the global convergence of the method
of section 4. In section 9 we formulate two more global convergence
theorems, which follow from results of sections 7 and 8. Numerical results
with simulated and experimental data are presented in Section 10 and 11
respectively. Concluding remarks are given in section 12.

\section{Statements of Forward and Inverse Problems}

\label{sec:2}

Below all functions are real valued ones. Let $b>1$ be a known number, $x\in 
\mathbb{R}$ be the spatial variable and the function $c(x)\in C^{3}(\mathbb{R%
})$, represents the spatially distributed dielectric constant. We assume
that 
\begin{align}
c(x)& \in \lbrack 1,b],x\in \mathbb{R},  \label{2.1} \\
c(x)=1& \mbox{ if }x\in (-\infty ,\varepsilon ]\cup \lbrack 1,\infty ),
\label{2.2}
\end{align}%
where $\varepsilon \in \left( 0,1\right) $ is a small number. Let $T$ be a
positive number. We consider the following Cauchy problem for a 1D
hyperbolic PDE with a variable coefficient in the principal part of the
operator:%
\begin{equation}
c(x)u_{tt}=u_{xx},\text{ }x\in \mathbb{R},t\in \left( 0,T\right) ,
\label{2.3}
\end{equation}%
\begin{equation}
u\left( x,0\right) =0,u_{t}\left( x,0\right) =\delta \left( x\right) .
\label{2.4}
\end{equation}%
The problem of finding the function $u(x,t)$ from conditions (\ref{2.3}), (%
\ref{2.4}) is our Forward Problem.

Let $\tau \left( x\right) $ be the travel time needed for the wave to travel
from the point source at $\left\{ 0\right\} $ to the point $x$, 
\begin{equation}
\tau \left( x\right) =\dint\limits_{0}^{x}\sqrt{c(s)}ds  \label{2.5}
\end{equation}%
By (\ref{2.5}) the following 1D analog of the eikonal equation is valid:%
\begin{equation}
\tau ^{\prime }\left( x\right) =\sqrt{c\left( x\right) }.  \label{2.50}
\end{equation}

Let $H(z),z\in \mathbb{R}$ be the Heaviside function,%
\begin{equation*}
H(z)=\left\{ 
\begin{array}{c}
1,\quad z>0, \\ 
0,\quad z<0.%
\end{array}%
\right.
\end{equation*}

\textbf{Lemma 2.1} \cite[Lemma 2.1]{KLNSN}. \emph{For }$x\geq 0$\emph{, the
function }$u(x,t)$\emph{\ has the form }%
\begin{equation}
u(x,t)=\frac{H(t-\tau \left( x\right) )}{2c^{1/4}(x)}+\hat{u}(x,t),
\label{2.6}
\end{equation}%
\emph{where the function }$\hat{u}\in C^{2}(t\geq \tau \left( x\right) ),%
\hat{u}(x,\tau \left( x\right) )=0.$\emph{\ Furthermore,}%
\begin{equation}
\lim_{t\rightarrow \tau ^{\div }\left( x\right) }u(x,t)=\frac{1}{2c^{1/4}(x)}%
.  \label{2.7}
\end{equation}

\textbf{Lemma 2.2} [Absorbing boundary condition] \cite{KLNSN,SKN1}. \emph{%
Let }$b>\varepsilon $\emph{\ be the number in (\ref{2.1}). Let }$x_{1}\geq b$%
\emph{\ and }$x_{2}\leq \varepsilon $\emph{\ be two arbitrary numbers. Then
the solution }$u(x,t)$\emph{\ of Forward Problem (\ref{2.3}), (\ref{2.4})
satisfies the absorbing boundary condition at }$x=x_{1}$\emph{\ and }$%
x=x_{2} $\emph{, i.e. }%
\begin{equation}
u_{x}(x_{1},t)+u_{t}(x_{1},t)=0\mbox{ for }t\in (0,T),  \label{2.21}
\end{equation}%
\begin{equation}
u_{x}(x_{2},t)-u_{t}(x_{2},t)=0\mbox{ for }t\in (0,T).  \label{2.22}
\end{equation}

We are interested in the following inverse problem:

\textbf{Coefficient Inverse Problem (CIP)}. \emph{Suppose that the following
two functions }$g_{0}(t)$\emph{, }$g_{1}(t)$\emph{\ are known: }%
\begin{equation}
u\left( \varepsilon ,t\right) =g_{0}\left( t\right) ,\quad u_{x}\left(
\varepsilon ,t\right) =g_{1}\left( t\right) ,\quad t\in \left( 0,T\right) .
\label{2.8}
\end{equation}%
\emph{Determine the function }$c(x)$\emph{\ for }$x\in (\varepsilon ,1)$%
\emph{, assuming that the number }$b>\varepsilon $\emph{\ in (\ref{2.1}) is
known.}

\textbf{Remark 2.1}. \emph{Note that only the function }$g_{0}\left(
t\right) $\emph{\ can be measured. As to the function }$g_{1}\left( t\right)
,$\emph{\ it follows from (\ref{2.22}) that }%
\begin{equation}
g_{1}\left( t\right) =g_{0}^{\prime }\left( t\right) .  \label{2.80}
\end{equation}%
\emph{We differentiate noisy function using the Tikhonov regularization
method \cite{T}. Since this method is well known, we do not describe it here.%
}

\section{A Boundary Value Problem for a Nonlinear PDE With Non-Local Terms}

\label{sec:3}

We now introduce a change of variable 
\begin{equation}
q(x,t)=u(x,t+\tau (x)).  \label{3.1}
\end{equation}%
We will consider the function $q(x,t)$ only for $t\geq 0.$ Using (\ref{2.50}%
) and (\ref{3.1}), we obtain 
\begin{equation}
q_{xx}-2q_{xt}\tau ^{\prime }-q_{t}\tau ^{\prime \prime }=0.  \label{3.3}
\end{equation}%
Furthermore, it follows from (\ref{2.7}) 
\begin{equation}
q(x,0)=\frac{1}{2c^{1/4}(x)}\neq 0.  \label{3.4}
\end{equation}

By (\ref{2.5}) and (\ref{3.4}) 
\begin{equation}
\tau ^{\prime \prime }(x)=-\frac{q_{x}(x,0)}{2q^{3}(x,0)}.  \label{3.5}
\end{equation}

Substituting (\ref{2.50}), (\ref{3.4}) and (\ref{3.5}) in (\ref{3.3}), we
obtain 
\begin{equation}
L(q)=q_{xx}-q_{xt}\frac{1}{2q^{2}(x,0)}+q_{t}\frac{q_{x}(x,0)}{2q^{3}(x,0)}%
=0.  \label{3.6}
\end{equation}%
Equation, (\ref{3.6}) is a nonlinear PDE with respect to the function $%
q(x,t) $ with nonlocal terms $q_{x}(x,0)$ and $q(x,0)$. We now need to
obtain boundary conditions for the function $q.$ By (\ref{2.2}) and (\ref%
{2.5}) $\tau \left( x\right) =x$ for $x\in \left[ 0,\varepsilon \right] .$
Hence, (\ref{2.8}), (\ref{2.80}) and (\ref{3.1}) lead to%
\begin{equation}
q\left( \varepsilon ,t\right) =g_{0}\left( t+\varepsilon \right)
,q_{x}\left( \varepsilon ,t\right) =2g_{0}^{\prime }\left( t+\varepsilon
\right) ,\text{ }t\in \left( 0,T\right) ,  \label{3.60}
\end{equation}

We will solve equation (\ref{3.6}) in the rectangle 
\begin{equation}
\Omega =\left\{ {(x,t):x\in (\varepsilon ,b),t\in (0,T)}\right\} .
\label{3.7}
\end{equation}%
By ({\ref{3.1})} 
\begin{equation}
q_{x}(x,t)=u_{x}(x,t+\tau (x))+u_{t}(x,t+\tau (x))\tau ^{\prime }(x)
\label{3.8}
\end{equation}%
By (\ref{2.1}), (\ref{2.2}) and (\ref{2.5}) $\tau ^{\prime }(b)=1$. Hence,
using (\ref{2.21}) and (\ref{3.8}), we obtain 
\begin{equation}
q_{x}(b,t)=0.  \label{3.9}
\end{equation}

It follows from (\ref{2.1}) and (\ref{3.4}) that 
\begin{equation}
\frac{1}{2b^{1/4}}\leq q(x,0)\leq \frac{1}{2},\quad x\in \lbrack \varepsilon
,b].  \label{3.10}
\end{equation}%
In addition, we need $q\in C^{2}(\overline{\Omega })$ and we also need to
bound the norm $\left\Vert q\right\Vert _{C^{2}(\overline{\Omega })}$ from
the above. Let $R>0$ be an arbitrary number. We define the set $B(R,g_{0})$
as 
\begin{equation}
B(R,g_{0})=\left\{ 
\begin{array}{c}
q\in H^{4}(\Omega ):\left\Vert q\right\Vert _{H^{4}(\Omega )}<R, \\ 
q(\varepsilon ,t+\varepsilon )=g_{0}(t+\varepsilon ),q_{x}(\varepsilon
,t)=2g_{0}^{\prime }\left( t+\varepsilon \right) , \\ 
q_{x}(b,t)=0, \\ 
\frac{1}{2b^{1/4}}\leq q\left( x,0\right) \leq 1/2,x\in \lbrack \varepsilon
,b].%
\end{array}%
\right.  \label{3.11}
\end{equation}%
We assume below that 
\begin{equation}
B(R,g_{0})\neq \varnothing .  \label{3.15}
\end{equation}%
By embedding theorem $\overline{B(R,g_{0})}\subset C^{2}(\overline{\Omega })$
and 
\begin{equation}
\left\Vert q\right\Vert _{C^{2}(\overline{\Omega })}\leq KR,\text{ }\forall
q\in \overline{B(R,g_{0})}.  \label{3.12}
\end{equation}%
where the constant $K=K\left( \Omega \right) >0$ depends only on the domain $%
\Omega .$ Using (\ref{3.10}), we define the function $q^{0}\left( x,0\right) 
$ as: 
\begin{eqnarray}
q^{0}\left( x,0\right) &=&\left\{ 
\begin{array}{c}
q\left( x,0\right) \text{ if }q\left( x,0\right) \in \left[ 1/\left(
2b^{1/4}\right) ,1/2\right] , \\ 
1/\left( 2b^{1/4}\right) ,\text{ if }q\left( x,0\right) <1/\left(
2b^{1/4}\right) , \\ 
1/2\text{ if }q\left( x,0\right) >1/2.%
\end{array}%
\right. \text{ }\forall q\in \overline{B(R,g_{0})},  \label{40.7} \\
\forall x &\in &\left[ \varepsilon ,b\right] .  \notag
\end{eqnarray}%
Then the function $q^{0}\left( x,0\right) $ is piecewise continuously
differentiable in $\left[ \varepsilon ,b\right] $ and by (\ref{3.12}) and (%
\ref{40.7}) 
\begin{equation}
\left\Vert q^{0}\left( x,0\right) \right\Vert _{C\left[ \varepsilon ,b\right]
},\max_{\left[ \varepsilon ,b\right] }\left\vert q_{x}^{0}\left( x,0\right)
\right\vert \leq KR,\text{ }\forall q\in B(R,g_{0}),  \label{40.8}
\end{equation}%
\begin{equation}
\frac{1}{2b^{1/4}}\leq q^{0}(x,0)\leq \frac{1}{2},\quad x\in \lbrack
\varepsilon ,b].  \label{3.13}
\end{equation}

Thus, (\ref{3.6}), (\ref{3.60}), (\ref{3.9}), (\ref{40.7}) and (\ref{3.13})
result in the following BVP for a nonlinear PDE with non-local terms:%
\begin{equation}
q_{xx}(x,t)-q_{xt}(x,t)\frac{1}{2\left( q^{0}(x,0)\right) ^{2}}+q_{t}(x,t)%
\frac{q_{x}(x,0)}{2\left( q^{0}(x,0)\right) ^{3}}=0\text{ in }\Omega ,
\label{3.14}
\end{equation}%
\begin{equation}
q(\varepsilon ,t)=g_{0}(t+\varepsilon ),\text{ }q_{x}(\varepsilon
,t)=2g_{0}^{\prime }(t+\varepsilon ),\text{ }q_{x}(b,t)=0.  \label{3.141}
\end{equation}%
Thus, we focus below on the numerical solution of the following problem:

\textbf{Problem 3.1.}\emph{\ Find a function }$q\in B(R,g_{0})$\emph{\
satisfying conditions (\ref{3.14}), (\ref{3.141}), where the function }$%
q^{0}(x,0)$ \emph{is defined in (\ref{40.7}).}

Suppose that we have solved Problem 3.1. Then, using (\ref{3.4}) and (\ref%
{40.7}), we set%
\begin{equation}
c\left( x\right) =\frac{1}{\left( 2q^{0}\left( x,0\right) \right) ^{4}}.
\label{3.16}
\end{equation}

\section{Numerical Method for Problem 3.1}

\label{sec:4}

\subsection{The function $q_{0}\left( x,t\right) $}

\label{sec:4.1}

We now find the first approximation $q_{0}(x,t)$ for the function $q(x,t).$
Using (\ref{2.2}), we choose $c\left( x\right) \equiv 1$ as the first guess
for the function $c\left( x\right) .$ Hence, by (\ref{3.4}), 
\begin{equation}
q_{0}(x,0)\equiv \frac{1}{2}.  \label{40.1}
\end{equation}%
We now need to find the function $q_{0}(x,t).$ To do this, drop the
nonlinear third term in the left hand side of equation (\ref{3.14}) and,
using (\ref{40.1}) and (\ref{40.7})\emph{, }set $1/\left( 2q^{0}(x,0)\right)
^{2}:=2$. Then (\ref{3.14}), (\ref{3.141})\emph{\ }become:%
\begin{equation}
q_{0xx}(x,t)-2q_{0xt}(x,t)=0\text{ in }\Omega ,  \label{400}
\end{equation}%
\begin{equation}
q_{0}(\varepsilon ,t)=g_{0}(t+\varepsilon ),\text{ }q_{0x}(\varepsilon
,t)=2g_{0}^{\prime }(t+\varepsilon ),q_{0x}(b,t)=0.  \label{401}
\end{equation}

BVP (\ref{400}), (\ref{401}) has overdetermined boundary conditions (\ref%
{401}). Typically, QRM\ works well for BVPs for PDEs with overdetermined
boundary conditions \cite{KlibAPNUM,KL}. Therefore, we solve BVP (\ref{400}%
), (\ref{401}) via CQRM. This means that we consider the following
minimization problem:

\textbf{Minimization Problem Number }$0$\textbf{.} \emph{Assuming (\ref{3.15}%
), minimize the functional }$J_{\lambda ,\beta }^{(0)}:\overline{B(R,g_{0})}%
\rightarrow \mathbb{R}$\emph{\ on the set }$,$\emph{\ }%
\begin{equation}
J_{\lambda ,\beta }^{(0)}(q_{0})=\dint\limits_{\Omega }\left(
q_{0xx}-2q_{0xt}\right) ^{2}e^{2\lambda \varphi }dxdt+\beta \Vert q_{0}\Vert
_{H^{4}(\Omega )}^{2},  \label{40.3}
\end{equation}%
\emph{where }$e^{2\lambda \varphi }$\emph{\ is the Carleman Weight Function
for the operator }$\partial _{x}-2\partial _{t}$ \cite{SKN1,SKN2} 
\begin{equation}
e^{2\lambda \varphi }=e^{-2\lambda \left( x+\alpha t\right) },  \label{40.4}
\end{equation}%
\emph{\ where }$\alpha \in \left( 0,1/2\right) $\emph{\ is the parameter,
and }$\beta \in \left( 0,1\right) $ \emph{is the regularization parameter.
Both parameters will be chosen later. }

Theorem 6.1 guarantees that, for appropriate values of parameters $\lambda
,\beta ,$ there exists unique minimizer $q_{0,\min }$\emph{\ }$\overline{%
B(R,g_{0})}$ of the functional $J_{\lambda ,\beta }^{(0)}(q_{0}).$

\subsection{The function $q_{n}\left( x,t\right) $ for $n\geq 1$}

\label{sec:4.2}

Assume that functionals $J_{\lambda ,\beta }^{(m)}(q_{m}):\overline{%
B(R,g_{0})}\rightarrow \mathbb{R}$ are defined for $m=0,...,n-1$ and their
minimizers functions $\left\{ q_{m,\min }\right\} _{m=0}^{n-1}\subset 
\overline{B(R,g_{0})}$ are constructed already, all for the same values of
parameters $\lambda ,\alpha ,\beta $. Replace in (\ref{3.14}) $q\left(
x,t\right) $ with $q_{n}\left( x,t\right) ,q^{0}(x,0)$ with $q_{n-1,\min
}^{0}(x,0),q_{x}\left( x,t\right) $ with $\partial _{x}q_{\left( n-1\right)
,\min }\left( x,t\right) $ and $q_{t}\left( x,t\right) $ with $\partial
_{t}q_{\left( n-1\right) ,\min }\left( x,t\right) .$ Then problem (\ref{3.14}%
), (\ref{3.141}) becomes a linear one with respect to the function $%
q_{n}\left( x,t\right) ,$%
\begin{equation}
q_{nxx}(x,t)-\frac{q_{nxt}(x,t)}{2\left( q_{n-1,\min }^{0}(x,0)\right) ^{2}}+%
\frac{\partial _{t}q_{\left( n-1\right) ,\min }\left( x,t\right) \partial
_{x}q_{\left( n-1\right) \min }(x,0)}{2\left( q_{n-1,\min }^{0}(x,0)\right)
^{3}}=0\text{ in }\Omega ,  \label{403}
\end{equation}%
\begin{equation}
q_{n}(\varepsilon ,t)=g_{0}(t+\varepsilon ),\text{ }q_{nx}(\varepsilon
,t)=2g_{0}^{\prime }(t+\varepsilon ),q_{nx}(b,t)=0.  \label{404}
\end{equation}

To solve problem (\ref{403}), (\ref{404}), we consider the following
minimization problem:

\textbf{Minimization Problem Number }$n$\textbf{.} \emph{Assuming (\ref{3.15}%
), minimize the functional }$J_{\lambda ,\beta }^{(n)}:H^{4}\left( \Omega
\right) \rightarrow \mathbb{R}$\emph{\ on the set }$B(R,g_{0}),$%
\begin{equation*}
J_{\lambda ,\beta }^{(n)}(q_{n})=\int_{\Omega }\Big(q_{nxx}(x,t)-\frac{%
q_{nxt}(x,t)}{2\big(q_{n-1,\min }^{0}(x,0)\big)^{2}}+\frac{\partial
_{t}q_{\left( n-1\right) ,\min }\left( x,t\right) \partial _{x}q_{\left(
n-1\right) ,\min }(x,0)}{2\big(q_{n-1,\min }^{0}(x,0)\big)^{3}}\Big)%
^{2}e^{2\lambda \varphi }dxdt
\end{equation*}%
\begin{equation}
+\beta \left\Vert q_{n}\right\Vert _{H^{4}\left( \Omega \right) }^{2}.
\label{405}
\end{equation}

Suppose that there exists unique minimizer $q_{n,\min }\in \overline{%
B(R,g_{0})}$ of the functional $J_{\lambda ,\beta }^{(n)}(q_{n}).$ Then,
following (\ref{3.16}), (\ref{40.7}) and (\ref{3.16}), we set%
\begin{equation}
c_{n}\left( x\right) =\frac{1}{\left( 2q_{n,\min }^{0}\left( x,0\right)
\right) ^{4}}.  \label{406}
\end{equation}%
The rest of the analytical part of this paper is devoted to the convergence
analysis of the iterative numerical method presented in this section.

\section{The Carleman Estimate}

\label{sec:5}

In this section we formulate the Carleman estimate, which is the main tool
of our construction. This estimate follows from Theorem 3.1 of \cite{KLNSN}
as well as from (\ref{40.8}). Let $q\left( x,t\right) \in \overline{%
B(R,g_{0})}$ be an arbitrary function and let the function $q^{0}\left(
x,0\right) $ be constructed from the function $q\left( x,t\right) $ as in (%
\ref{40.7}). Consider the operator $L_{0}:H^{2}(\Omega )\rightarrow
L_{2}(\Omega ),$%
\begin{equation*}
L_{0}u=u_{xx}(x,t)-u_{xt}(x,t)\frac{1}{2\left( q^{0}(x,0)\right) ^{2}},\quad %
\mbox{for all }(x,t)\in \Omega .
\end{equation*}

\textbf{Theorem 5.1} (Carleman estimate \cite{KLNSN}). \emph{There exists a
number }$\alpha _{0}=\alpha _{0}\left( R,\Omega \right) >0$\emph{\ depending
only on }$R,\Omega $ \emph{such that for any }$\alpha \in (0,\alpha _{0})$%
\emph{\ there exists a sufficiently large number }$\lambda _{0}=\lambda
_{0}\left( R,\Omega ,\alpha \right) >1$\emph{\ depending only on }$R$,$%
\Omega ,\alpha $\emph{\ such that for all }$\lambda \geq \lambda _{0}$\emph{%
\ and for all functions }$v\in H^{2}(\Omega )$\emph{\ the following Carleman
estimate holds:}%
\begin{equation*}
\dint\limits_{\Omega }\left( L_{0}v\right) ^{2}e^{2\lambda \varphi }dxdt
\end{equation*}%
\begin{equation*}
\geq C\dint\limits_{\Omega }\big(\lambda \left( v_{x}^{2}+v_{t}^{2}\right)
+\lambda ^{3}v^{2}\big)e^{2\lambda \varphi }dxdt+C\dint\limits_{\varepsilon
}^{b}\big(\lambda v_{x}^{2}(x,0)+\lambda ^{3}v^{2}(x,0)\big)e^{-2\lambda x}dx
\end{equation*}%
\begin{equation}
-C\dint\limits_{0}^{T}\left( \lambda \left( v_{x}^{2}(\varepsilon
,t)+v_{t}^{2}(\varepsilon ,t)\right) +\lambda ^{3}v^{2}(\varepsilon
,t)\right) e^{-2\lambda (\varepsilon +\alpha t)}dt  \label{4.3}
\end{equation}%
\begin{equation*}
-C\dint\limits_{\varepsilon }^{b}\left( \lambda v_{x}^{2}(x,T)+\lambda
^{3}v^{2}(x,T)\right) e^{-2\lambda (x+\alpha T)}dx.
\end{equation*}%
\emph{\ }

\textbf{Remark 5.1.} \emph{Here and everywhere below }$C=C\left( R,\Omega
,\alpha \right) >0$\emph{\ denotes different constants depending only on
listed parameters.}

\section{Strict Convexity of Functional (\protect\ref{405}) on the Set $%
B(R,g_{0})$, Existence and Uniqueness of Its Minimizer}

\label{sec:6}

Functional (\ref{405}) is quadratic. We prove in this section that it is
strictly convex on the set $\overline{B(R,g_{0})}.$ In addition, we prove
existence and uniqueness of its minimizer on this set. Although similar
results were proven in many of the above cited publications on the
convexification, see, e.g. \cite{KLNSN} for the closest one, there are some
peculiarities here, which are important for our convergence analysis, see
Remarks 6.1, 6.2 below.

Introduce the subspace $H_{0}^{4}\left( \Omega \right) \subset H^{4}\left(
\Omega \right) $ as:%
\begin{equation}
H_{0}^{4}\left( \Omega \right) =\left\{ v\in H^{4}\left( \Omega \right)
:v\left( \varepsilon ,t\right) =v_{x}\left( \varepsilon ,t\right)
=v_{x}\left( b,t\right) =0\right\} .  \label{5.1}
\end{equation}%
Denote $\left[ ,\right] $ the scalar product in the space $H^{4}\left(
\Omega \right) .$ Also, denote 
\begin{equation}
A_{n}\left( q\right) \left( x,t\right) =q_{xx}(x,t)-q_{xt}(x,t)\frac{1}{%
2\left( q_{n-1,\min }^{0}(x,0)\right) ^{2}},  \label{5.100}
\end{equation}%
\begin{equation}
B_{n}\left( x,t\right) =\partial _{t}q_{\left( n-1\right) ,\min }\left(
x,t\right) \frac{\partial _{x}q_{\left( n-1\right) ,\min }(x,0)}{2\left(
q_{n-1,\min }^{0}(x,0)\right) ^{3}}.  \label{5.101}
\end{equation}

\textbf{Theorem 6.1.}\emph{\ 1. For any set of parameters }$\lambda ,\beta $%
\emph{\ and } \emph{for any }$q\in \overline{B(R,g_{0})}$ \emph{the
functional }$J_{\lambda ,\beta }^{\left( n\right) }$\emph{\ has the
Fr\'echet derivative }$\left( J_{\lambda ,\beta }^{\left( n\right)
}(q)\right) ^{\prime }\in H_{0}^{4}\left( \Omega \right) .$\emph{\ The
formula for }$\left( J_{\lambda ,\beta }^{\left( n\right) }(q)\right)
^{\prime }$ is:%
\begin{equation}
\left( J_{\lambda ,\beta }^{\left( n\right) }(q)\right) ^{\prime }\left(
h\right) =2\dint\limits_{\Omega }\left( A_{n}\left( q\right) \left(
x,t\right) +B_{n}\left( x,t\right) \right) A_{n}\left( h\right) \left(
x,t\right) e^{2\lambda \varphi }dxdt +2\beta \left[ q,h\right],  \label{5.2}
\end{equation}
for all $h\in H_{0}^{4}\left( \Omega \right).$

2. \emph{This derivative is Lipschitz continuous in }$\overline{B(R,g_{0})}$%
\emph{, i.e. there exists a constant }$D>0$\emph{\ such that }%
\begin{equation}
\left\Vert \left( J_{\lambda ,\beta }^{\left( n\right) }(q_{2})\right)
^{\prime }-\left( J_{\lambda ,\beta }^{\left( n\right) }(q_{1})\right)
^{\prime }\right\Vert _{H^{4}\left( \Omega \right) }\leq D\left\Vert
q_{2}-q_{1}\right\Vert _{H^{4}\left( \Omega \right) },\quad \mbox{for all }%
q_{1},q_{2}\in \overline{B(R,g_{0})}.  \label{5.3}
\end{equation}

\emph{3. Let }$\alpha _{0}=\alpha _{0}\left( R,\Omega \right) >0,\alpha \in
\left( 0,\alpha _{0}\right) $\emph{\ and }$\lambda _{0}=\lambda _{0}\left(
R,\Omega ,\alpha \right) \geq 1$ \emph{be the numbers of Theorem 5.1}$.$%
\emph{\ Then there exists a sufficiently large constant }%
\begin{equation}
\lambda _{1}=\lambda _{1}\left( R,\Omega ,\alpha \right) \geq \lambda _{0}>1
\label{5.4}
\end{equation}%
\emph{depending only on listed parameters such that for all }$\lambda \geq
\lambda _{1}$ \emph{and for all\ }$\beta \in \left[ 2e^{-\lambda \alpha
T},1\right) $\emph{\ the functional }$J_{\lambda ,\beta }^{\left( n\right)
}(q)$\emph{\ is strictly convex on the set }$\overline{B\left( R,g\right) }.$%
\emph{\ More precisely, let }$q\in \overline{B(R,g)}$ \emph{be an arbitrary
function and also let the function }$q+h\in \overline{B(R,g)}.$\emph{\ Then} 
\emph{the following inequality holds:}%
\begin{equation*}
J_{\lambda ,\beta }^{\left( n\right) }(q+h)-J_{\lambda ,\beta }^{\left(
n\right) }(q)-\left( J_{\lambda ,\beta }^{\left( n\right) }(q)\right)
^{\prime }(h)\geq C\dint\limits_{\Omega }\Big[\lambda \left(
h_{x}^{2}+h_{t}^{2}\right) +\lambda ^{3}h^{2}\Big]e^{2\lambda \varphi }dxdt
\end{equation*}%
\begin{equation}
+C\dint\limits_{\varepsilon }^{b}\big(\lambda h_{x}^{2}(x,0)+\lambda
^{3}h^{2}(x,0)\big)e^{-2\lambda x}dx+\frac{\beta }{2}\left\Vert h\right\Vert
_{H^{4}\left( \Omega \right) }^{2},\text{ }\forall \lambda \geq \lambda _{1}.
\label{5.5}
\end{equation}

\emph{4. For any }$\lambda \geq \lambda _{1}$ \emph{there exists unique
minimizer }%
\begin{equation}
q_{n,\min }\in \overline{B(R,g_{0})}  \label{5.50}
\end{equation}%
\emph{\ of the functional }$J_{\lambda ,\beta }^{\left( n\right) }(q)$\emph{%
\ on the set }$\overline{B(R,g_{0})}$\emph{\ and the following inequality
holds:}%
\begin{equation}
\left[ \left( J_{\lambda ,\beta }^{\left( n\right) }(q)\right) ^{\prime
},q-q_{n,\min }\right] \geq 0,\text{ }\forall q\in \overline{B(R,g_{0})},
\label{5.6}
\end{equation}%
\emph{where }$\left[ ,\right] $\emph{\ denotes the scalar product in }$%
H^{4}\left( \Omega \right) .$

\textbf{Remark 6.1}. \emph{Even though the expression in the right hand side
of (\ref{5.2}) is linear with respect to the function }$q,$\emph{\ we cannot
use Riesz theorem here to prove existence and uniqueness of the minimizer }$%
q_{n,\min },$\emph{\ at which }$\left( J_{\lambda ,\beta }^{\left( n\right)
}(q_{n,\min })\right) ^{\prime }=0.$\emph{\ Rather, all what we can prove is
(\ref{5.6}). This is because we need to ensure that the function }$q_{n,\min
}\in \overline{B(R,g_{0})}$\emph{. }

\textbf{Proof of Theorem 6.1}. Since both function $q,q+h\in \overline{B(R,g)%
}$ satisfy the same boundary conditions, then 
\begin{equation}
h\in H_{0}^{4}\left( \Omega \right) .  \label{5.7}
\end{equation}%
By (\ref{405}) and (\ref{5.7}) 
\begin{equation*}
J_{\lambda ,\beta }^{(n)}(q+h)-J_{\lambda ,\beta
}^{(n)}(q)=2\dint\limits_{\Omega }\left( A_{n}\left( q\right) \left(
x,t\right) +B_{n}\left( x,t\right) \right) A_{n}\left( h\right) \left(
x,t\right) e^{2\lambda \varphi }dxdt+2\beta \left[ q,h\right]
\end{equation*}%
\begin{equation}
+\dint\limits_{\Omega }\left[ A_{n}\left( h\right) \left( x,t\right) \right]
^{2}e^{2\lambda \varphi }dxdt+\beta \left\Vert h\right\Vert _{H^{4}\left(
\Omega \right) }^{2},\text{ }\forall h\in H_{0}^{4}\left( \Omega \right) .
\label{5.8}
\end{equation}%
The expression in the first line of (\ref{5.8}) coincides with expression in
the right hand side of (\ref{5.2}). In fact, this is a bounded linear
functional mapping $H_{0}^{4}\left( \Omega \right) $ in $\mathbb{R}$.
Therefore, by Riesz theorem, there exists unique function $\widetilde{J}%
_{n}\left( q\right) \in H_{0}^{4}\left( \Omega \right) $ such that%
\begin{equation*}
\left[ \widetilde{J}_{n}\left( q\right) ,h\right] =2\dint\limits_{\Omega
}\left( A_{n}\left( q\right) \left( x,t\right) +B_{n}\left( x,t\right)
\right) A_{n}\left( h\right) \left( x,t\right) e^{2\lambda \varphi
}dxdt+2\beta \left[ q,h\right] ,\text{ }\forall h\in H_{0}^{4}\left( \Omega
\right) .
\end{equation*}%
In addition, it is clear from (\ref{5.8}) and (\ref{5.9}) that%
\begin{equation*}
\lim_{\left\Vert h\right\Vert _{H^{4}\left( \Omega \right) }\rightarrow
0}\left\{ \frac{1}{\left\Vert h\right\Vert _{H^{4}\left( \Omega \right) }}%
\left[ J_{\lambda ,\beta }^{(n)}(q+h)-J_{\lambda ,\beta }^{(n)}(q)-\left[ 
\widetilde{J}_{n}\left( q\right) ,h\right] \right] \right\} =0.
\end{equation*}%
Therefore, 
\begin{equation}
\widetilde{J}_{n}\left( q\right) =\left( J_{\lambda ,\beta }^{\left(
n\right) }(q)\right) ^{\prime }\in H_{0}^{4}\left( \Omega \right)
\label{5.10}
\end{equation}%
is the Fr\'echet derivative of the functional $J_{\lambda ,\beta }^{\left(
n\right) }(q):\overline{B(R,g_{0})}\rightarrow \mathbb{R}$ at the point $%
q\in \overline{B(R,g_{0})},$ and the right hand side of (\ref{5.2}) indeed
represents $\left( J_{\lambda ,\beta }^{\left( n\right) }(q)\right) ^{\prime
}\left( h\right) .$ Estimate (\ref{5.3}) obviously follows from (\ref{5.2}).

We now prove strict convexity estimate (\ref{5.5}). To do this, we apply
Carleman estimate (\ref{4.3}) to the third line of (\ref{5.8}). We obtain%
\begin{multline}
J_{\lambda ,\beta }^{(n)}(q+h)-J_{\lambda ,\beta }^{(n)}(q)-\left(
J_{\lambda ,\beta }^{\left( n\right) }(q)\right) ^{\prime }\left( h\right) \\
\geq C\dint\limits_{\Omega }\big(\lambda \left( h_{x}^{2}+h_{t}^{2}\right)
+\lambda ^{3}h^{2}\big)e^{2\lambda \varphi }dxdt+C\dint\limits_{\varepsilon
}^{b}\big(\lambda h_{x}^{2}(x,0)+\lambda ^{3}h^{2}(x,0)\big)e^{-2\lambda x}dx
\\
-C\dint\limits_{\varepsilon }^{b}\left( \lambda h_{x}^{2}(x,T)+\lambda
^{3}h^{2}(x,T)\right) e^{-2\lambda (x+\alpha T)}dx+\beta \left\Vert
h\right\Vert _{H^{4}\left( \Omega \right) }^{2},\text{ }\forall h\in
H_{0}^{4}\left( \Omega \right) .  \label{5.11}
\end{multline}%
By trace theorem there exists a constant $C_{1}=C_{1}\left( \Omega \right)
>0 $ depending only on the domain $\Omega $ such that 
\begin{equation*}
\left\Vert v\right\Vert _{H^{4}\left( \Omega \right) }^{2}\geq
C_{1}\left\Vert v\left( x,T\right) \right\Vert _{H^{1}\left( \Omega \right)
}^{2},\text{ }\forall v\in H^{4}\left( \Omega \right) .
\end{equation*}%
Since the regularization parameter $\beta \in \left[ 2e^{-\lambda \alpha
T},1\right) ,$ then we can choose $\lambda _{1}$ so large that 
\begin{equation*}
\frac{\beta }{2}C_{1}\geq C_{1}e^{-\lambda \alpha T}>Ce^{-2\lambda (x+\alpha
T)},\text{ }\forall \lambda \geq \lambda _{1},\forall x\in \left[
\varepsilon ,b\right] .
\end{equation*}%
Hence, for these values of $\lambda ,$ the expression in the last line of (%
\ref{5.11}) can be estimated from the below as:%
\begin{equation}
-C\dint\limits_{\varepsilon }^{b}\left( \lambda h_{x}^{2}(x,T)+\lambda
^{3}h^{2}(x,T)\right) e^{-2\lambda (x+\alpha T)}dx+\beta \left\Vert
h\right\Vert _{H^{4}\left( \Omega \right) }^{2}\geq \frac{\beta }{2}%
\left\Vert h\right\Vert _{H^{4}\left( \Omega \right) }^{2},\text{ }\forall
h\in H_{0}^{4}\left( \Omega \right) .  \label{5.12}
\end{equation}%
Hence, (\ref{5.11}) and (\ref{5.12}) imply%
\begin{multline*}
J_{\lambda ,\beta }^{(n)}(q+h)-J_{\lambda ,\beta }^{(n)}(q)-\left(
J_{\lambda ,\beta }^{\left( n\right) }(q)\right) ^{\prime }\left( h\right)
\geq C\dint\limits_{\Omega }\big(\lambda \left( h_{x}^{2}+h_{t}^{2}\right)
+\lambda ^{3}h^{2}\big)e^{2\lambda \varphi }dxdt \\
+C\dint\limits_{\varepsilon }^{b}\big(\lambda h_{x}^{2}(x,0)+\lambda
^{3}h^{2}(x,0)\big)e^{-2\lambda x}dx+\frac{\beta }{2}\left\Vert h\right\Vert
_{H^{4}\left( \Omega \right) }^{2},\forall \lambda \geq \lambda _{1}.
\end{multline*}%
This proves (\ref{5.5}). The existence and uniqueness of the minimizer $%
q_{\min ,n}\in \overline{B(R,g_{0})}$ and inequality (\ref{5.6})\emph{\ }%
follow from (\ref{5.5}) as well as from a combination of Lemma 2.1 and
Theorem 2.1 of \cite{BAK}, also see \cite[Chapter 10, section 3]{Minoux}. $%
\square $

\textbf{Remark 6.2}. \emph{Since the functional }$J_{\lambda ,\beta
}^{\left( n\right) }(q)$\emph{\ is quadratic, then its strict convexity on
the whole space }$H^{4}\left( \Omega \right) $\emph{\ follows immediately
from the presence of the regularization term }$\beta \left\Vert
q_{n}\right\Vert _{H^{4}\left( \Omega \right) }^{2}$\emph{\ in it. However,
in addition to the claim of its strict convexity, we actually need in our
convergence analysis below those terms in the right hand side of the strict
convexity estimate (\ref{5.5}), which are different from the term }$\beta
\left\Vert h\right\Vert _{H^{4}\left( \Omega \right) }^{2}/2.$\emph{\ These
terms are provided by Carleman estimate (\ref{5.5}). The condition }$\beta
\in \left[ 2e^{-\lambda \alpha T},1\right) $\emph{\ of Theorem 6.1 is
imposed to dominate the negative term in the last line of (\ref{5.5}).}

\section{How to Find the Minimizer}

\label{sec:7}

Since we search for the minimizer $q_{n,\min }$ of functional (\ref{405}) on
the bounded set $\overline{B(R,g_{0})}$ rather than on the whole space $%
H^{4}\left( \Omega \right) ,$ then we cannot just use Riesz theorem to find
this minimizer, also see \cite{BAK} and \cite[Chapter 10, section 3]{Minoux}
for the case of finding a minimizer of a strictly convex functional on a
bounded set. Two ways of finding the minimizer $q_{n,\min }\in \overline{%
B(R,g_{0})}$ are described in this section.

\subsection{Gradient descent method}

\label{sec:7.1}

Keeping in mind (\ref{3.15}), choose an arbitrary function $q_{0,n}\in
B(R,g_{0}).$ We arrange the gradient descent method of the minimization of
functional (\ref{405}) as follows:%
\begin{equation}
q_{k,n}=q_{\left( k-1\right) ,n}-\eta \left( J_{\lambda ,\beta }^{\left(
n\right) }(q_{\left( k-1\right) ,n})\right) ^{\prime },\text{ }k=1,2,...,
\label{7.1}
\end{equation}%
where $\eta \in \left( 0,1\right) $ is a small number, which is chosen
later. It is important to note that since functions $J_{\lambda ,\beta
}^{\left( n\right) }(q_{\left( k-1\right) ,n})\in H_{0}^{4}\left( \Omega
\right) ,$ then boundary conditions (\ref{404}) are kept the same for all
functions $q_{k,n}\left( x,t\right) .$ Also, using (\ref{40.7}) and (\ref%
{406}), we set%
\begin{equation}
c_{k,n}\left( x\right) =\frac{1}{\left( 2q_{k,n}^{0}\left( x,0\right)
\right) ^{4}},\text{ }x\in \left[ \varepsilon ,b\right] ,  \label{7.2}
\end{equation}%
\begin{equation}
c_{n,\min }\left( x\right) =\frac{1}{\left( 2q_{n,\min }^{0}\left(
x,0\right) \right) ^{4}},\text{ }x\in \left[ \varepsilon ,b\right] .
\label{7.3}
\end{equation}

Theorem 7.1 claims the global convergence of the gradient descent method (%
\ref{7.1}) in the case when $q_{\min ,n}\in B(R/3,g_{0}).$

\textbf{Theorem 7.1}. \emph{Let the number }$\lambda _{1}=\lambda _{1}\left(
R,\Omega ,\alpha \right) \geq \lambda _{0}>1$\emph{\ be the one defined in (%
\ref{5.4}). Let }$\lambda \geq \lambda _{1}$\emph{. For this value of }$%
\lambda ,$\emph{\ let }$q_{n,\min }\in \overline{B\left( R,g_{0}\right) }$%
\emph{\ be the unique minimizer of the functional }$J_{\lambda ,\beta
}^{\left( n\right) }(q_{n})$\emph{\ on the set }$\overline{B\left(
R,g_{0}\right) }$\emph{\ with the regularization parameter }$\beta \in \left[
2e^{-\lambda \alpha T},1\right) $\emph{\ (Theorem 6.1)}$.$ \emph{Assume that
the function }$q_{n,\min }\in B(R/3,g_{0}).$\emph{\ For each }$n,$ \emph{%
choose the starting point of the gradient descent method (\ref{7.1}) as }$%
q_{0,n}\in B(R/3,g_{0}).$ \emph{Then, there exists a number }$\eta _{0}\in
\left( 0,1\right) $\emph{\ such that for any }$\eta \in (0,\eta _{0})$\emph{%
\ all functions }$q_{k,n}\in B(R,g_{0}).$\emph{\ Furthermore, there exists a
number }$\theta =\theta \left( \eta \right) \in (0,1)$\emph{\ such that the
following convergence estimates are valid: }%
\begin{equation}
\left\Vert q_{k,n}-q_{n,\min }\right\Vert _{H^{4}\left( \Omega \right) }\leq
\theta ^{k}\left\Vert q_{0,n}-q_{n,\min }\right\Vert _{H^{4}\left( \Omega
\right) },k=1,...,  \label{7.4}
\end{equation}%
\begin{equation}
\left\Vert c_{k,n}-c_{n,\min }\right\Vert _{H^{3}\left( \varepsilon
,b\right) }\leq C\theta ^{k}\left\Vert q_{0,n}-q_{n,\min }\right\Vert
_{H^{4}\left( \Omega \right) },\text{ }k=1,...  \label{7.5}
\end{equation}

\textbf{Proof}. Estimate (\ref{7.4}) follows immediately from \cite[Theorem
4.6]{KlibSAR2} combined with \textquotedblleft corrections" of functions $%
q_{k,n}\left( x,0\right) ,q_{n,\min }\left( x,0\right) $ via (\ref{40.7}).
Estimate (\ref{7.5}) follows immediately from trace theorem, (\ref{40.7})
and (\ref{7.2})-(\ref{7.4}). $\square $

\subsection{Gradient projection method}

\label{sec:7.2}

Suppose now that there is no information on whether or not the function 
\emph{\ }$q_{n,\min }\in B(R/3,g_{0}).$ In this case we construct the
gradient projection method. We introduce the function $F\left( x,t\right) $
since it is easy to construct the projection operator on a ball with the
center at $\left\{ 0\right\} .$

Suppose that there exists a function $F\in H^{4}\left( \Omega \right) $ such
that%
\begin{equation}
F(\varepsilon ,t)=g_{0}(t+\varepsilon ),\text{ }F_{x}(\varepsilon
,t)=2g_{0}^{\prime }(t+\varepsilon ),\text{ }F_{x}(b,t)=0.  \label{7.50}
\end{equation}%
This function exists, for example, in the case when the function $%
g_{0}(t)\in H^{5}\left( 0,T+\varepsilon \right) .$ Indeed, as one of
examples of this function (among many), one can take, e.g. 
\begin{equation*}
F\left( x,t\right) =\chi \left( x\right) \left( g_{0}(t+\varepsilon
)+2xg_{0}^{\prime }(t+\varepsilon )\right) ,
\end{equation*}%
where the function $\chi \left( x\right) $ is such that 
\begin{equation*}
\chi \left( x\right) \in C^{4}\left[ \varepsilon ,b\right] ,\chi \left(
x\right) =\left\{ 
\begin{array}{c}
1,x\in \left[ \varepsilon ,b/4\right] , \\ 
0,x\in \left[ b/2,b\right] , \\ 
\text{between }0\text{ and }1\text{ for }x\in \left( b/4,b/2\right) .%
\end{array}%
\right.
\end{equation*}%
The existence of such functions $\chi \left( x\right) $ is well known from
the Analysis course. Denote%
\begin{equation}
p_{n}\left( x,t\right) =q_{n}\left( x,t\right) -F(x,t).  \label{7.6}
\end{equation}%
Then (\ref{403}), (\ref{404}) become:%
\begin{equation}
p_{nxx}(x,t)-p_{nxt}(x,t)\frac{1}{2\left( q_{n-1,\min }^{0}(x,0)\right) ^{2}}
\label{7.7}
\end{equation}%
\begin{equation*}
+F_{xx}(x,t)-F_{xt}(x,t)\frac{1}{2\left( q_{n-1,\min }^{0}(x,0)\right) ^{2}}%
+\partial _{t}q_{\left( n-1\right) ,\min }\left( x,t\right) \frac{\partial
_{x}q_{\left( n-1\right) ,\min }(x,0)}{2\left( q_{n-1,\min }^{0}(x,0)\right)
^{3}}=0\text{ in }\Omega ,
\end{equation*}%
\begin{equation}
p_{n}(\varepsilon ,t)=0,\text{ }p_{nx}(\varepsilon ,t)=0,p_{nx}(b,t)=0.
\label{7.8}
\end{equation}%
Assume that 
\begin{equation}
\left\Vert F\right\Vert _{H^{4}\left( \Omega \right) }<R.  \label{7.80}
\end{equation}%
By (\ref{7.6}), (\ref{7.8}), (\ref{7.80}) and triangle inequality%
\begin{equation}
p_{n}\in B_{0}\left( 2R\right) =\left\{ p\in H_{0}^{4}\left( \Omega \right)
:\left\Vert p\right\Vert _{H_{0}^{4}\left( \Omega \right) }<2R\right\} .
\label{7.81}
\end{equation}%
To find the function $p_{n}\in B_{0}\left( 2R\right) $ satisfying conditions
(\ref{7.7}), (\ref{7.8}), we minimize the following functional $I_{\lambda
,\beta }^{\left( n\right) }(p_{n}):\overline{B_{0}\left( 2R\right) }%
\rightarrow \mathbb{R}$%
\begin{equation}
I_{\lambda ,\beta }^{\left( n\right) }(p_{n})=J_{\lambda ,\beta }^{\left(
n\right) }(p_{n}+F),\text{ }p_{n}\in B_{0}\left( 2R\right) .  \label{7.9}
\end{equation}

\textbf{Remark 7.1.} \emph{An obvious analog of Theorem 6.1 is valid of
course for functional (\ref{7.9}). But in this case we should have instead
of (\ref{5.4}) }$\lambda \geq \widetilde{\lambda }_{1}=\lambda _{1}\left(
2R,\Omega ,\alpha \right) \geq \lambda _{0}>1.$ \emph{In particular, there
exists unique minimizer }$p_{n,\min }\in \overline{B_{0}\left( 2R\right) }$%
\emph{\ of this functional on the closed ball }$\overline{B_{0}\left(
2R\right) }.$\emph{\ We omit the formulation of this theorem, since it is an
obvious reformulation of Theorem 6.1.}

Let $P_{\overline{B_{0}\left( 2R\right) }}:H_{0}^{4}\left( \Omega \right)
\rightarrow \overline{B_{0}\left( 2R\right) }$ be the projection operator of
the space $H_{0}^{4}\left( \Omega \right) $ on the closed ball $\overline{%
B_{0}\left( 2R\right) }\subset H_{0}^{4}\left( \Omega \right) .$ Then this
operator can be easily constructed:%
\begin{equation*}
P_{\overline{B_{0}\left( 2R\right) }}\left( p\right) =\left\{ 
\begin{array}{c}
p\text{ if }p\in \overline{B_{0}\left( 2R\right) }, \\ 
\frac{2R}{\left\Vert p\right\Vert _{H^{4}\left( \Omega \right) }}p\text{ if }%
p\notin \overline{B_{0}\left( 2R\right) }.%
\end{array}%
\right.
\end{equation*}

We now construct the gradient projection method of the minimization of the
functional $I_{\lambda ,\beta }^{\left( n\right) }(p_{n})$ on the set $%
\overline{B_{0}\left( 2R\right) }.$ Let $p_{0,n}\in B_{0}\left( 2R\right) $
be an arbitrary function. Then the sequence of the gradient projection
method is:%
\begin{equation}
p_{k,n}=P_{\overline{B_{0}\left( 2R\right) }}\left( p_{\left( k-1\right)
,n}-\eta \left( I_{\lambda ,\beta }^{\left( n\right) }(p_{\left( k-1\right)
,n})\right) ^{\prime }\right) ,\text{ }k=1,2,...,  \label{7.10}
\end{equation}%
where $\eta \in \left( 0,1\right) $ is a small number, which is chosen
later. Using (\ref{7.2}), (\ref{7.3}) and (\ref{7.6}), we set%
\begin{equation}
\widetilde{c}_{k,n}\left( x\right) =\frac{1}{\left( 2\widetilde{q}%
_{k,n}^{0}\left( x,0\right) \right) ^{4}},\text{ }x\in \left[ \varepsilon ,b%
\right] ,  \label{7}
\end{equation}%
\begin{equation}
\widetilde{c}_{n,\min }\left( x\right) =\frac{1}{\left( 2\widetilde{q}%
_{n,\min }^{0}\left( x,0\right) \right) ^{4}},\text{ }x\in \left[
\varepsilon ,b\right] .  \label{8}
\end{equation}%
Here the function $\widetilde{q}_{n,k}^{0}\left( x,0\right) $ is obtained as
follows: First, we consider the function $\left( p_{n,k}+F\right) \left(
x,0\right) .$ Next, we apply procedure (\ref{40.7}) to this function.
Similarly for $\widetilde{q}_{n,\min }^{0}\left( x,0\right) .$

Denote $p_{n,\min }\in \overline{B_{0}\left( 2R\right) }$ the unique
minimizer of functional (\ref{7.9}) on the set $\overline{B_{0}\left(
2R\right) }$ (Remark 7.1). Following (\ref{7.6}), denote 
\begin{equation}
\widetilde{q}_{n,\min }=p_{n,\min }+F,\text{ }\widetilde{q}_{k,n}=p_{k,n}+F.
\label{9}
\end{equation}%
We omit the proof of Theorem 7.2 since it is very similar with the proof of
Theorem 7.1. The only difference is that instead of Theorem 4.6 of \cite%
{KlibSAR2} one should use Theorem 4.1 of \cite{KLNSN}.

\textbf{Theorem 7.2}. \emph{Let (\ref{7.50}), (\ref{7.80}) hold. Let the
number }$\lambda _{1}=\lambda _{1}\left( R,\Omega ,\alpha \right) \geq
\lambda _{0}>1$\emph{\ be the one defined in (\ref{5.4}) and let }%
\begin{equation*}
\lambda \geq \widetilde{\lambda }_{1}=\lambda _{1}\left( 2R,\Omega ,\alpha
\right) \geq \lambda _{1}\left( R,\Omega ,\alpha \right) .
\end{equation*}%
\emph{For this value of }$\lambda $ \emph{and for }$\beta \in \left[
2e^{-\lambda \alpha T},1\right) $\emph{\ let }$p_{n,\min }\in \overline{%
B_{0}\left( 2R\right) }$\emph{\ be the unique minimizer of functional (\ref%
{7.9}) on the set }$\overline{B_{0}\left( 2R\right) }$ \emph{(Remark 7.1).
Let notations (\ref{7}) hold. For each }$n,$ \emph{choose the starting point
of the gradient projection method (\ref{7.1}) as }$p_{0,n}\in B_{0}(2R).$ 
\emph{Then there exists a number }$\eta _{0}\in \left( 0,1\right) $\emph{\
such that for any }$\eta \in (0,\eta _{0})$\emph{\ there exists a number }$%
\theta =\theta \left( \eta \right) \in (0,1)$\emph{\ such that the following
convergence estimates are valid for the iterative process (\ref{7.10}): }%
\begin{equation*}
\left\Vert \widetilde{q}_{k,n}-\widetilde{q}_{n,\min }\right\Vert
_{H^{4}\left( \Omega \right) }\leq \theta ^{k}\left\Vert \widetilde{q}_{0,n}-%
\widetilde{q}_{n,\min }\right\Vert _{H^{4}\left( \Omega \right) },
\end{equation*}%
\begin{equation*}
\left\Vert \widetilde{c}_{k,n}-\widetilde{c}_{n,\min }\right\Vert
_{H^{4}\left( \Omega \right) }\leq C\theta ^{k}\left\Vert \widetilde{q}%
_{0,n}-\widetilde{q}_{n,\min }\right\Vert _{H^{4}\left( \Omega \right) },
\end{equation*}%
\emph{where functions }$\widetilde{c}_{k,n},\widetilde{c}_{n,\min },%
\widetilde{q}_{k,n},\widetilde{q}_{n,\min }$\emph{\ are defined in (\ref{7}%
)-(\ref{9}).}

\textbf{Remark 7.2. }\emph{Both theorems 7.1 and 7.2 claim the global
convergence of corresponding minimization procedures. This is because the
starting function in both cases is an arbitrary one either in }$B\left(
R/3,g_{0}\right) $\emph{\ or in }$B_{0}\left( 2R\right) $\emph{\ and a
smallness condition is not imposed on the number }$R$\emph{.}

\section{Contraction Mapping and Global Convergence}

\label{sec:8}

In this section, we prove the global convergence of the numerical method of
section 4 for solving Problem 3.1. To do this, we introduce first the exact
solution of our CIP. Recall that the concept of the existence of the exact
solution is one of the main concepts of the theory of ill-posed problems 
\cite{BK,T}. In particular, an estimate in our global convergence theorem is
very similar to the one in contraction mapping.

Suppose that there exists a function $c^{\ast }\left( x\right) $ satisfying
conditions (\ref{2.1}), (\ref{2.2}). Let $u^{\ast }\left( x,t\right) $ be
the solution of problem (\ref{2.3}), (\ref{2.4}) with $c:=c^{\ast }.$ We
assume that the corresponding data $g_{0}^{\ast }\left( t\right) ,\partial
_{t}g_{0}^{\ast }\left( t\right) $ for the CIP are noiseless, see (\ref{2.8}%
), (\ref{2.80}). Let $q^{\ast }\left( x,t\right) $ be the function $q^{\ast
}\left( x,t\right) $ which is constructed from the function $u^{\ast }\left(
x,t\right) $ as in (\ref{3.1}). Following (\ref{3.11}), we assume that 
\begin{equation}
q^{\ast }\in B(R,g_{0}^{\ast }),  \label{8.1}
\end{equation}%
\begin{equation}
B(R,g_{0}^{\ast })=\left\{ 
\begin{array}{c}
q\in H^{4}(\Omega ):\left\Vert q\right\Vert _{H^{4}(\Omega )}<R, \\ 
q(\varepsilon ,t+\varepsilon )=g_{0}^{\ast }(t+\varepsilon
),q_{x}(\varepsilon ,t)=2\left( g_{0}^{\ast }\right) ^{\prime }\left(
t+\varepsilon \right) , \\ 
q_{x}(b,t)=0, \\ 
\frac{1}{2b^{1/4}}\leq q\left( x,0\right) \leq 1/2,x\in \lbrack \varepsilon
,b]%
\end{array}%
\right\} .  \label{8.2}
\end{equation}%
By (\ref{3.14})-(\ref{3.141}) 
\begin{equation}
q_{xx}^{\ast }(x,t)-q_{xt}^{\ast }(x,t)\frac{1}{2\left( q^{\ast
}(x,0)\right) ^{2}}+q_{t}^{\ast }(x,t)\frac{q_{x}^{\ast }(x,0)}{2\left(
q^{\ast }(x,0)\right) ^{3}}=0\text{ in }\Omega ,  \label{8.3}
\end{equation}%
\begin{equation}
q^{\ast }(\varepsilon ,t)=g_{0}^{\ast }(t+\varepsilon ),\text{ }q_{x}^{\ast
}(\varepsilon ,t)=2\partial _{t}g_{0}^{\ast }\left( t+\varepsilon \right)
,q_{x}^{\ast }(b,t)=0.  \label{8.4}
\end{equation}%
By (\ref{3.4})%
\begin{equation}
c^{\ast }\left( x\right) =\frac{1}{\left( 2q^{\ast }\left( x,0\right)
\right) ^{4}}.  \label{8.40}
\end{equation}

It is important in the formulation of Theorem 6.1 that both functions $q$
and $q+h$ should have the same boundary conditions as prescribed in $%
B(R,g_{0}).$ However, boundary conditions for functions $q_{n}$ and $q^{\ast
}$ are different. Hence, similarly to subsection 7.2, we consider a function 
$F^{\ast }\in H^{4}\left( \Omega \right) $ such that (see (\ref{7.50}))%
\begin{equation}
F^{\ast }(\varepsilon ,t)=g_{0}^{\ast }(t+\varepsilon ),\text{ }F_{x}^{\ast
}(\varepsilon ,t)=2\partial _{t}g_{0}^{\ast }(t+\varepsilon ),\text{ }%
F_{x}^{\ast }(b,t)=0.  \label{8.6}
\end{equation}%
We assume similarly to (\ref{7.80}) that%
\begin{equation}
\left\Vert F^{\ast }\right\Vert _{H^{4}\left( \Omega \right) }<R.
\label{8.7}
\end{equation}%
Also, similarly to (\ref{7.6}), we introduce the function $p^{\ast }\left(
x,t\right) $ as: 
\begin{equation}
p^{\ast }\left( x,t\right) =q^{\ast }\left( x,t\right) -F^{\ast }(x,t).
\label{8.8}
\end{equation}%
Let $\lambda _{1}$ be the number defined in (\ref{5.4}), let $\lambda \geq
\lambda _{1}$ and the function $q_{n,\min }\in \overline{B\left(
R,g_{0}\right) }$ (see (\ref{5.50})) be the unique minimizer of the
functional $J_{\lambda ,\beta }^{\left( n\right) }(q_{n})$ on the set $%
\overline{B\left( R,g_{0}\right) },$ the existence of which is established
in Theorem 6.1. Following (\ref{7.6}), denote%
\begin{equation}
p_{n,\min }\left( x,t\right) =q_{n,\min }\left( x,t\right) -F\left(
x,t\right) .  \label{8.9}
\end{equation}%
By (\ref{5.50}), (\ref{7.81}), (\ref{8.1}), (\ref{8.2}) and (\ref{8.7})-(\ref%
{8.9})%
\begin{equation}
p_{n,\min },p^{\ast }\in \overline{B_{0}\left( 2R\right) }.  \label{8.10}
\end{equation}%
Also, it follows from embedding theorem, (\ref{3.12}), (\ref{5.50}), (\ref%
{7.80}), (\ref{8.7}) and (\ref{8.10}) that 
\begin{equation}
\left\Vert q^{\ast }\right\Vert _{C^{2}\left( \overline{\Omega }\right)
},\left\Vert p^{\ast }\right\Vert _{C^{2}\left( \overline{\Omega }\right)
},\left\Vert q_{n,\min }\right\Vert _{_{C^{2}\left( \overline{\Omega }%
\right) }},\left\Vert p_{n,\min }\right\Vert _{_{C^{2}\left( \overline{%
\Omega }\right) }}\leq C.  \label{8.11}
\end{equation}

We assume that the data $g_{0},g_{0}^{\prime }$ for our CIP are given with
noise of the level $\delta ,$ where the number $\delta >0$ is sufficiently
small. More precisely, we assume that 
\begin{equation}
\left\Vert F-F^{\ast }\right\Vert _{H^{4}\left( \Omega \right) }<\delta .
\label{8.12}
\end{equation}

Observe that (\ref{40.7}), (\ref{8.1}) and (\ref{8.2}) imply that%
\begin{equation}
\left\vert q_{n-1,\min }^{0}\left( x,0\right) -q^{\ast }\left( x,0\right)
\right\vert \leq \left\vert q_{n-1,\min }\left( x,0\right) -q^{\ast }\left(
x,0\right) \right\vert ,\text{ }x\in \left[ \varepsilon ,b\right] .
\label{8.13}
\end{equation}%
By (\ref{8.3}), (\ref{8.4}), (\ref{8.6}) and (\ref{8.8})%
\begin{equation}
p_{xx}^{\ast }(x,t)-p_{xt}^{\ast }(x,t)\frac{1}{2\left( q^{\ast
}(x,0)\right) ^{2}}+q_{t}^{\ast }(x,t)\frac{q_{x}^{\ast }(x,0)}{2\left(
q^{\ast }(x,0)\right) ^{3}}\text{ }+F_{xx}^{\ast }-F_{xt}^{\ast }(x,t)\frac{1%
}{2\left( q^{\ast }(x,0)\right) ^{2}}=0  \label{8.14}
\end{equation}%
in $\Omega \times \lbrack 0,T]$ and 
\begin{equation}
p^{\ast }(\varepsilon ,t)=p_{x}^{\ast }(\varepsilon ,t)=p_{x}^{\ast }(b,t)=0.
\label{8.15}
\end{equation}

\textbf{Theorem 8.1} (contraction mapping and global convergence of the
method of section 3). \emph{Let functions }$F$\emph{,}$F^{\ast }\in
H^{4}\left( \Omega \right) $\emph{\ satisfy conditions (\ref{7.50}), (\ref%
{7.80}), (\ref{8.6}), (\ref{8.7}) and (\ref{8.12}). Let a sufficiently large
number }$\lambda _{1}=\lambda _{1}\left( R,\Omega ,\alpha \right) \geq
\lambda _{0}>1$\emph{\ be the one defined in (\ref{5.4}). Let }%
\begin{equation}
\lambda \geq \widetilde{\lambda }_{1}=\lambda _{1}\left( 2R,\Omega ,\alpha
\right) \geq \lambda _{1}\left( R,\Omega ,\alpha \right) .  \label{8.151}
\end{equation}%
\emph{For this value of }$\lambda ,$\emph{\ let }$q_{n,\min }\in \overline{%
B\left( R,g_{0}\right) }$\emph{\ be the unique minimizer of the functional }$%
J_{\lambda ,\beta }^{\left( n\right) }(q_{n})$\emph{\ on the set }$\overline{%
B\left( R,g_{0}\right) }$\emph{\ with the regularization parameter }$\beta
\in \left[ 2e^{-\lambda \alpha T},1\right) $\emph{\ (Theorem 6.1)}$.$\emph{\
Let }%
\begin{equation}
\overline{q}_{n}=q_{n,\min }-q^{\ast },\text{ }\overline{c}_{n}=c_{n,\min
}-c^{\ast },\emph{\ }  \label{10}
\end{equation}%
\emph{where }$c_{n,\min }$\emph{\ is defined in (\ref{7.3}). Then the
following convergence estimate holds}%
\begin{equation*}
\dint\limits_{\Omega }\left( \overline{q}_{nx}^{2}+\overline{q}_{nt}^{2}+%
\overline{q}_{n}^{2}\right) \left( x,t\right) e^{2\lambda \varphi
}dxdt+\dint\limits_{\varepsilon }^{b}\left( \overline{q}_{nx}^{2}+\overline{q%
}_{n}^{2}\right) \left( x,0\right) e^{-2\lambda x}dx
\end{equation*}%
\begin{equation}
\leq \frac{C}{\lambda }\dint\limits_{\Omega }\left( \overline{q}_{\left(
n-1\right) x}^{2}+\overline{q}_{\left( n-1\right) t}^{2}+\overline{q}%
_{n-1}^{2}\right) \left( x,t\right) e^{2\lambda \varphi }dxdt  \label{8.16}
\end{equation}%
\begin{equation*}
+\frac{C}{\lambda }\dint\limits_{\varepsilon }^{b}\left( \overline{q}%
_{\left( n-1\right) x}^{2}\left( x,0\right) +\overline{q}_{n-1}^{2}\right)
\left( x,0\right) e^{-2\lambda x}dx+C\left( \delta ^{2}+\beta \right) ,
\end{equation*}%
\emph{which leads to:}%
\begin{equation}
\dint\limits_{\Omega }\left( \overline{q}_{nx}^{2}+\overline{q}_{nt}^{2}+%
\overline{q}_{n}^{2}\right) \left( x,t\right) e^{2\lambda \varphi }dxdt
\label{8.17}
\end{equation}%
\begin{equation*}
\leq \frac{C^{n}}{\lambda ^{n}}\dint\limits_{\Omega }\left( \overline{q}%
_{0x}^{2}+\overline{q}_{0t}^{2}+\overline{q}_{0}^{2}\right) \left(
x,t\right) e^{2\lambda \varphi }dxdt+C\left( \delta ^{2}+\beta \right) .
\end{equation*}%
\emph{In addition, \ }%
\begin{equation}
\left\Vert \overline{c}_{n}\right\Vert _{H^{1}\left( \varepsilon ,b\right)
}^{2}\leq \frac{C^{n}}{\lambda ^{n}}\dint\limits_{\Omega }\left( \overline{q}%
_{0x}^{2}+\overline{q}_{0t}^{2}+\overline{q}_{0}^{2}\right) \left(
x,t\right) e^{2\lambda \varphi }dxdt+C\left( \delta ^{2}+\beta \right) .
\label{8.18}
\end{equation}

\textbf{Remarks 8.1}:

1. \emph{Due to the presence of the term }$C/\lambda $\emph{\ with a
sufficiently large }$\lambda $\emph{, estimate (\ref{8.16}) is similar to
the one in the classic contraction mapping principle, although we do not
claim here the existence of the fixed point. This explains the title of our
paper. }

\emph{2. When computing the unique minimizer }$q_{0,\min }$\emph{\ of
functional (\ref{40.3}) on the set }$\overline{B(R,g_{0})},$\emph{\ we do
not impose a smallness condition on the number }$R$\emph{. Therefore,
Theorem 8.1 claims the global convergence of the method of section 4, see
Introduction for our definition of the global convergence. The same is true
for Theorems 9.1 and 9.2 in the next section.}

\textbf{Proof Theorem 8.1. }Denote $h_{n}=p^{\ast }-p_{n,\min }.$ By (\ref%
{8.9}) $h_{n}=-\overline{q}_{n}+\left( F-F^{\ast }\right) .$ Hence, (\ref%
{8.12}) and embedding theorem imply: 
\begin{equation}
h_{n}^{2}\geq \frac{1}{2}\overline{q}_{n}^{2}-C\delta ^{2},h_{nx}^{2}\geq 
\frac{1}{2}\overline{q}_{nx}^{2}-C\delta ^{2},h_{nt}^{2}\geq \frac{1}{2}%
\overline{q}_{nt}^{2}-C\delta ^{2}\text{ in }\overline{\Omega },
\label{8.190}
\end{equation}%
\begin{equation}
h_{n}^{2}+h_{nx}^{2}+h_{nt}^{2}\leq C\left( \overline{q}_{n}^{2}+\overline{q}%
_{nx}^{2}+\overline{q}_{nt}^{2}+\delta ^{2}\right) \text{ in }\overline{%
\Omega }.  \label{8.191}
\end{equation}

Consider the functional $I_{\lambda ,\beta }^{\left( n\right) }(p^{\ast
})=J_{\lambda ,\beta }^{\left( n\right) }(p^{\ast }+F).$ Since both
functions $p_{n,\min }$ and $p^{\ast }$ have the same zero boundary
conditions (\ref{7.8}) and since by (\ref{8.10}) both of them belong to the
set $\overline{B_{0}\left( 2R\right) }$, then the analog of Theorem 6.1,
which is mentioned in Remark 7.1, implies (see (\ref{5.5}))%
\begin{equation*}
I_{\lambda ,\beta }^{\left( n\right) }(p^{\ast })-I_{\lambda ,\beta
}^{\left( n\right) }(p_{n,\min })-\left( I_{\lambda ,\beta }^{\left(
n\right) }(p_{n,\min })\right) ^{\prime }(h_{n})
\end{equation*}%
\begin{equation}
\geq C\dint\limits_{\Omega }\Big[\lambda \left( h_{nx}^{2}+h_{nt}^{2}\right)
+\lambda ^{3}h_{n}^{2}\Big]e^{2\lambda \varphi }dxdt  \label{8.20}
\end{equation}%
\begin{equation*}
+C\dint\limits_{\varepsilon }^{b}\big(\lambda h_{nx}^{2}(x,0)+\lambda
^{3}h_{n}^{2}(x,0)\big)e^{-2\lambda x}dx+\frac{\beta }{2}\left\Vert
h_{n}\right\Vert _{H^{4}\left( \Omega \right) }^{2},\text{ }\forall \lambda
\geq \widetilde{\lambda }_{1}.
\end{equation*}%
By (\ref{5.6}) 
\begin{equation*}
-\left( I_{\lambda ,\beta }^{\left( n\right) }(p_{n,\min })\right) ^{\prime
}(h_{n})\leq 0.
\end{equation*}%
Hence, the left hand side of (\ref{8.20}) can be estimated as:%
\begin{equation}
I_{\lambda ,\beta }^{\left( n\right) }(p^{\ast })-I_{\lambda ,\beta
}^{\left( n\right) }(p_{n,\min })-\left( I_{\lambda ,\beta }^{\left(
n\right) }(p_{n,\min })\right) ^{\prime }(h)\leq I_{\lambda ,\beta }^{\left(
n\right) }(p^{\ast }).  \label{8.21}
\end{equation}

We now estimate the right hand side of (\ref{8.21}) from the above. It
follows from (\ref{405}), (\ref{7.9}) and (\ref{8.3}) that 
\begin{equation}
I_{\lambda ,\beta }^{(n)}(p^{\ast })=\dint\limits_{\Omega
}G_{n}^{2}e^{2\lambda \varphi }dxdt+\beta \left\Vert p^{\ast }+F\right\Vert
_{H^{4}\left( \Omega \right) }^{2},  \label{8.22}
\end{equation}%
where 
\begin{equation*}
G_{n}=p_{xx}^{\ast }(x,t)-p_{xt}^{\ast }(x,t)\frac{1}{2\left( q_{\left(
n-1\right) ,\min }^{0}(x,0)\right) ^{2}}+\partial _{t}q_{\left( n-1\right)
,\min }\left( x,t\right) \frac{\partial _{x}q_{\left( n-1\right) ,\min }(x,0)%
}{2\left( q_{\left( n-1\right) ,\min }^{0}(x,0)\right) ^{3}}
\end{equation*}%
\begin{equation*}
+F_{xx}-F_{xt}(x,t)\frac{1}{2\left( q_{\left( n-1\right) ,\min
}^{0}(x,0)\right) ^{2}}
\end{equation*}%
\begin{equation}
=q_{xx}^{\ast }(x,t)-q_{xt}^{\ast }(x,t)\frac{1}{2\left( q^{\ast
}(x,0)\right) ^{2}}+q_{t}^{\ast }(x,t)\frac{q_{x}^{\ast }(x,0)}{2\left(
q^{\ast }(x,0)\right) ^{3}}  \label{8.23}
\end{equation}%
\begin{equation*}
+\left( F-F^{\ast }\right) _{xx}-\left( F-F^{\ast }\right) _{xt}\frac{1}{%
2\left( q_{\left( n-1\right) ,\min }^{0}(x,0)\right) ^{2}}
\end{equation*}%
\begin{equation*}
-q_{xt}^{\ast }\left( \frac{1}{2\left( q_{\left( n-1\right) ,\min
}^{0}(x,0)\right) ^{2}}-\frac{1}{2\left( q^{\ast }(x,0)\right) ^{2}}\right)
\end{equation*}%
\begin{equation*}
+\left( \partial _{t}q_{\left( n-1\right) ,\min }\left( x,t\right) \frac{%
\partial _{x}q_{\left( n-1\right) ,\min }(x,0)}{2\left( q_{\left( n-1\right)
,\min }^{0}(x,0)\right) ^{3}}-q_{t}^{\ast }(x,t)\frac{q_{x}^{\ast }(x,0)}{%
2\left( q^{\ast }(x,0)\right) ^{3}}\right) .
\end{equation*}%
By (\ref{8.3}), the third line of (\ref{8.23}) equals zero. Hence, (\ref%
{8.23}) becomes 
\begin{equation*}
G_{n}=\left( F-F^{\ast }\right) _{xx}+\left( F-F^{\ast }\right) _{xt}\frac{1%
}{2\left( q_{\left( n-1\right) ,\min }^{0}(x,0)\right) ^{2}}
\end{equation*}%
\begin{equation*}
-q_{xt}^{\ast }\left( \frac{1}{2\left( q_{\left( n-1\right) ,\min
}^{0}(x,0)\right) ^{2}}-\frac{1}{2\left( q^{\ast }(x,0)\right) ^{2}}\right)
\end{equation*}%
\begin{equation*}
+\left( \partial _{t}q_{\left( n-1\right) ,\min }\left( x,t\right) \frac{%
\partial _{x}q_{\left( n-1\right) ,\min }(x,0)}{2\left( q_{\left( n-1\right)
,\min }^{0}(x,0)\right) ^{3}}-q_{t}^{\ast }(x,t)\frac{q_{x}^{\ast }(x,0)}{%
2\left( q^{\ast }(x,0)\right) ^{3}}\right) .
\end{equation*}%
Hence, by (\ref{40.7}), (\ref{8.12}) and embedding theorem%
\begin{equation*}
\left\vert G_{n}\right\vert \leq C\delta +C\left\vert \frac{1}{2\left(
q_{\left( n-1\right) ,\min }^{0}(x,0)\right) ^{2}}-\frac{1}{2\left( q^{\ast
}(x,0)\right) ^{2}}\right\vert
\end{equation*}%
\begin{equation}
+\left\vert \partial _{t}q_{\left( n-1\right) ,\min }\left( x,t\right) \frac{%
\partial _{x}q_{\left( n-1\right) ,\min }(x,0)}{2\left( q_{\left( n-1\right)
,\min }^{0}(x,0)\right) ^{3}}-q_{t}^{\ast }(x,t)\frac{q_{x}^{\ast }(x,0)}{%
2\left( q^{\ast }(x,0)\right) ^{3}}\right\vert .  \label{8.24}
\end{equation}%
Using (\ref{8.12}) and (\ref{8.13}), we obtain 
\begin{equation*}
\left\vert \frac{1}{2\left( q_{\left( n-1\right) ,\min }^{0}(x,0)\right) ^{2}%
}-\frac{1}{2\left( q^{\ast }(x,0)\right) ^{2}}\right\vert
\end{equation*}%
\begin{equation*}
=\frac{\left\vert \left( q_{\left( n-1\right) ,\min }^{0}\left( x,0\right)
-q^{\ast }(x,0)\right) -\left( F-F^{\ast }\right) \left( x,0\right)
\right\vert \left\vert q_{\left( n-1\right) ,\min }^{0}(x,0)+q^{\ast
}(x,0)\right\vert }{2\left( q_{\left( n-1\right) ,\min }^{0}(x,0)\right)
^{2}\left( q^{\ast }(x,0)\right) ^{2}.}
\end{equation*}%
\begin{equation*}
\leq C\delta +C\left\vert h_{n-1}\left( x,0\right) \right\vert .
\end{equation*}%
Combining this with (\ref{8.24}), we obtain%
\begin{equation*}
\left\vert G_{n}\right\vert \leq C\delta +C\left\vert h_{n-1}\left(
x,0\right) \right\vert
\end{equation*}%
\begin{equation}
+\left\vert \partial _{t}q_{\left( n-1\right) ,\min }\left( x,t\right) \frac{%
\partial _{x}q_{\left( n-1\right) ,\min }(x,0)}{2\left( q_{\left( n-1\right)
,\min }^{0}(x,0)\right) ^{3}}-q_{t}^{\ast }(x,t)\frac{q_{x}^{\ast }(x,0)}{%
2\left( q^{\ast }(x,0)\right) ^{3}}\right\vert .  \label{8.25}
\end{equation}%
Next, 
\begin{equation*}
\frac{1}{2\left( q_{\left( n-1\right) ,\min }^{0}(x,0)\right) ^{3}}=\frac{1}{%
2\left( q^{\ast }(x,0)\right) ^{3}}
\end{equation*}%
\begin{equation*}
+\left( \frac{1}{2\left( q_{\left( n-1\right) ,\min }^{0}\left( x,0\right)
\right) ^{3}}-\frac{1}{2\left( q^{\ast }(x,0)\right) ^{3}}\right)
\end{equation*}%
\begin{equation*}
=\frac{1}{2\left( q^{\ast }(x,0)\right) ^{3}}+S_{n-1}\left( x,0\right) \left[
\left( q_{\left( n-1\right) ,\min }^{0}\left( x,0\right) -q^{\ast
}(x,0)\right) -\left( F-F^{\ast }\right) \left( x,0\right) \right] ,
\end{equation*}%
where the function $S_{n-1}\left( x,0\right) $ can be estimated as%
\begin{equation}
\left\vert S_{n-1}\left( x,0\right) \right\vert \leq C.  \label{8.26}
\end{equation}%
Hence, 
\begin{equation*}
\partial _{t}q_{\left( n-1\right) ,\min }\left( x,t\right) \frac{\partial
_{x}q_{\left( n-1\right) ,\min }(x,0)}{2\left( q_{\left( n-1\right) ,\min
}^{0}(x,0)\right) ^{3}}=
\end{equation*}%
\begin{equation*}
\frac{\partial _{t}q_{\left( n-1\right) ,\min }\left( x,t\right) \partial
_{x}q_{\left( n-1\right) ,\min }(x,0)}{2\left( q^{\ast }(x,0)\right) ^{3}}+
\end{equation*}%
\begin{equation*}
+\left[ \partial _{t}q_{\left( n-1\right) ,\min }\left( x,t\right) \partial
_{x}q_{\left( n-1\right) ,\min }(x,0)\right] \times
\end{equation*}%
\begin{equation*}
\times S_{n-1}\left( x,0\right) \left[ \left( q_{\left( n-1\right) ,\min
}^{0}\left( x,0\right) -q^{\ast }(x,0)\right) -\left( F-F^{\ast }\right)
\left( x,0\right) \right] .
\end{equation*}%
Hence, using (\ref{3.11}), (\ref{40.7}), (\ref{8.1}), (\ref{8.13}) and (\ref%
{8.26}), we obtain 
\begin{equation*}
\left\vert \partial _{t}q_{\left( n-1\right) ,\min }\left( x,t\right) \frac{%
\partial _{x}q_{\left( n-1\right) ,\min }(x,0)}{2\left( q_{\left( n-1\right)
,\min }^{0}(x,0)\right) ^{3}}-q_{t}^{\ast }(x,t)\frac{q_{x}^{\ast }(x,0)}{%
2\left( q^{\ast }(x,0)\right) ^{3}}\right\vert
\end{equation*}%
\begin{equation}
\leq \frac{1}{2\left( q^{\ast }(x,0)\right) ^{3}}\left\vert \partial
_{t}q_{\left( n-1\right) ,\min }\left( x,t\right) \partial _{x}q_{\left(
n-1\right) ,\min }(x,0)-q_{t}^{\ast }(x,t)q_{x}^{\ast }(x,0)\right\vert
\label{8.27}
\end{equation}%
\begin{equation*}
+C\left( \delta +\left\vert h_{n-1}\left( x,0\right) \right\vert \right)
\end{equation*}%
\begin{equation*}
\leq C\left( \delta +\left\vert h_{n-1}\left( x,0\right) \right\vert
+\left\vert \partial _{x}h_{n-1,x}\left( x,0\right) \right\vert +\left\vert
h_{n-1,t}\left( x,t\right) \right\vert \right) .
\end{equation*}%
Combining (\ref{8.25}) with (\ref{8.27}), we obtain 
\begin{equation*}
\left\vert G_{n}\left( x,t\right) \right\vert \leq C\left( \delta
+\left\vert h_{n-1}\left( x,0\right) \right\vert +\left\vert h_{n-1,x}\left(
x,0\right) \right\vert +\left\vert h_{n-1,t}\left( x,t\right) \right\vert
\right) .
\end{equation*}%
Hence, by (\ref{8.22}) 
\begin{equation*}
I_{\lambda ,\beta }^{(n)}(p^{\ast })\leq C\dint\limits_{\Omega }\left(
\delta ^{2}+h_{n-1}^{2}\left( x,0\right) +h_{n-1,x}^{2}\left( x,0\right)
+h_{n-1,t}^{2}\left( x,t\right) \right) e^{2\lambda \varphi }dxdt+C\beta .
\end{equation*}%
Substituting this in (\ref{8.21}) and then using (\ref{8.20}), we obtain%
\begin{equation*}
\dint\limits_{\Omega }\left( h_{nx}^{2}+h_{nt}^{2}+h_{n}^{2}\right)
e^{2\lambda \varphi }dxdt+\dint\limits_{\varepsilon }^{b}\left(
h_{nx}^{2}(x,0)+h_{n}^{2}(x,0)\right) e^{-2\lambda x}dx
\end{equation*}%
\begin{equation}
\leq \frac{C}{\lambda }\dint\limits_{\Omega }\left( \delta
^{2}+h_{n-1}^{2}\left( x,0\right) +h_{n-1,x}^{2}\left( x,0\right)
+h_{n-1,t}^{2}\left( x,t\right) \right) e^{2\lambda \varphi }dxdt+C\beta .
\label{8.28}
\end{equation}%
Obviously%
\begin{equation}
\dint\limits_{\Omega }\left( h_{n-1}^{2}\left( x,0\right)
+h_{n-1,x}^{2}\left( x,0\right) \right) e^{2\lambda \varphi }dxdt\leq \frac{1%
}{2\lambda \alpha }\dint\limits_{\varepsilon }^{b}\left( h_{\left(
n-1\right) x}^{2}(x,0)+h_{n-1}^{2}(x,0)\right) e^{-2\lambda x}dx,
\label{8.29}
\end{equation}%
\begin{equation}
\dint\limits_{\Omega }\delta ^{2}e^{2\lambda \varphi }dxdt\leq C\frac{\delta
^{2}}{\lambda ^{2}}.  \label{8.30}
\end{equation}%
Denote 
\begin{equation}
y_{n}=\dint\limits_{\Omega }\left( h_{nx}^{2}+h_{nt}^{2}+h_{n}^{2}\right)
e^{2\lambda \varphi }dxdt+\dint\limits_{\varepsilon }^{b}\left(
h_{nx}^{2}(x,0)+h_{n}^{2}(x,0)\right) e^{-2\lambda x}dx.  \label{8.31}
\end{equation}%
Then (\ref{8.28})-(\ref{8.30}) imply%
\begin{equation}
y_{n}\leq \frac{C}{\lambda }y_{n-1}+C\left( \frac{\delta ^{2}}{\lambda ^{2}}%
+\beta \right) .  \label{8.32}
\end{equation}%
Iterating (\ref{8.32}) with respect to $n,$ we obtain%
\begin{equation*}
\dint\limits_{\Omega }\left( h_{nx}^{2}+h_{nt}^{2}+h_{n}^{2}\right)
e^{2\lambda \varphi }dxdt+\dint\limits_{\varepsilon }^{b}\left(
h_{nx}^{2}(x,0)+h_{n}^{2}(x,0)\right) e^{-2\lambda x}dx
\end{equation*}%
\begin{equation}
\leq \frac{C}{\lambda ^{n}}\left[ \dint\limits_{\Omega }\left(
h_{0x}^{2}+h_{0t}^{2}+h_{0}^{2}\right) e^{2\lambda \varphi
}dxdt+\dint\limits_{\varepsilon }^{b}\left(
h_{0x}^{2}(x,0)+h_{0}^{2}(x,0)\right) e^{-2\lambda x}dx\right]  \label{8.33}
\end{equation}%
\begin{equation*}
+C\left( \frac{\delta ^{2}}{\lambda ^{2}}+\beta \right) .
\end{equation*}%
Apply (\ref{8.190}) to the left hand side of estimate (\ref{8.33}). Also,
apply (\ref{8.191}) at $n=0$ to the right hand side of (\ref{8.33}). We
obtain (\ref{8.17}). Estimate (\ref{8.18}) follows from an obvious
combination of (\ref{8.17}) with (\ref{7.3}), (\ref{8.40}), (\ref{8.13}) and
(\ref{8.150}). Finally, estimate (\ref{8.16}) follows immediately from (\ref%
{8.190}), (\ref{8.191}) and (\ref{8.28}). $\square $

\section{Global Convergence of the Gradient and Gradient Projection Methods
to the Exact Solution}

\label{sec:9}

First, we consider the gradient method of the minimization of functionals $%
J_{\lambda ,\beta }^{\left( n\right) }(q_{\left( k-1\right) ,n})$ on the set 
$\overline{B(R,g_{0})},$ see (\ref{7.1}). The proof of Theorem 9.1 follows
immediately from the triangle inequality combined with Theorems 7.1 and 8.1.

\textbf{Theorem 9.1}. \emph{Let }$\alpha _{0}$\emph{\ and }$\lambda _{0}$%
\emph{\ be the numbers of Theorem 5.1.} \emph{Let the sufficiently large
number }$\lambda _{1}=\lambda _{1}\left( R,\Omega ,\alpha \right) \geq
\lambda _{0}>1$\emph{\ be the one defined in (\ref{5.4}). Let the number }$%
\widetilde{\lambda }_{1}$\emph{\ be the same as in (\ref{8.151}), }%
\begin{equation*}
\widetilde{\lambda }_{1}=\lambda _{1}\left( 2R,\Omega ,\alpha \right) \geq
\lambda _{1}\left( R,\Omega ,\alpha \right) .
\end{equation*}%
\emph{Let }$\lambda \geq \lambda _{1}$\emph{\ and let the regularization
parameter }$\beta \in \left[ 2e^{-\lambda \alpha T},1\right) .$\emph{\
Assume that the functions }$q_{\min ,n}\in B(R/3,g_{0})$\emph{\ for all }$n.$%
\emph{\ For each} $n$, \emph{choose the starting point of the gradient
method (\ref{7.1}) as }$q_{0,n}\in B(R/3,g_{0}).$ \emph{Then there exists a
number }$\eta _{0}\in \left( 0,1\right) $\emph{\ such that for any }$\eta
\in (0,\eta _{0})$\emph{\ all functions }$q_{k,n}\in B(R,g_{0}),k=1,...;$ $%
n=1,....$\emph{\ Furthermore, there exists a number }$\theta =\theta \left(
\eta \right) \in (0,1)$\emph{\ such that the following convergence estimate
is valid:}%
\begin{equation*}
\left\Vert c_{k,n}-c^{\ast }\right\Vert _{H^{1}\left( \varepsilon ,b\right)
}\leq C\theta ^{k}\left\Vert q_{0,n}-q_{n,\min }\right\Vert _{H^{4}\left(
\Omega \right) }
\end{equation*}%
\begin{equation*}
+\frac{C^{n/2}}{\lambda ^{n/2}}\left( \dint\limits_{\Omega }\left( \overline{%
q}_{0x}^{2}+\overline{q}_{0t}^{2}+\overline{q}_{0}^{2}\right) \left(
x,t\right) e^{2\lambda \varphi }dxdt\right) ^{1/2}+C\left( \delta +\sqrt{%
\beta }\right) ,
\end{equation*}%
\emph{where the function }$\overline{q}_{0}$\emph{\ is defined in (\ref{10}).%
}

Consider now the gradient projection method of the minimization of the
functionals $I_{\lambda ,\beta }^{\left( n\right) }(p_{n})=J_{\lambda ,\beta
}^{\left( n\right) }(p_{n}+F)$ in (\ref{7.9}) on the set $\overline{%
B_{0}\left( 2R\right) }$, see (\ref{7.10}).\emph{\ }We use notations (\ref{7}%
). Theorem 9.2 follows immediately from the triangle inequality combined
with Theorems 7.2 and 8.1.

\textbf{Theorem 9.2.} \emph{Let the number }$\lambda _{1}=\lambda _{1}\left(
R,\Omega ,\alpha \right) \geq \lambda _{0}>1$\emph{\ be the one defined in (%
\ref{5.4}). Let the number }$\widetilde{\lambda }_{1}$\emph{\ be the same as
in (\ref{8.151}), }%
\begin{equation*}
\widetilde{\lambda }_{1}=\lambda _{1}\left( 2R,\Omega ,\alpha \right) \geq
\lambda _{1}\left( R,\Omega ,\alpha \right) .
\end{equation*}%
\emph{Let }$\lambda \geq \lambda _{1}$\emph{\ and let the regularization
parameter }$\beta \in \left[ 2e^{-\lambda \alpha T},1\right) .$\emph{\
Consider the gradient projection method (\ref{7.10}). For each} $n$, \emph{%
choose the starting point }$p_{0,n}$ \emph{of this method as an arbitrary
point of the ball }$B_{0}(2R).$ \emph{Then there exists a number }$\eta
_{0}\in \left( 0,1\right) $\emph{\ such that for any }$\eta \in (0,\eta
_{0}) $\emph{\ there exists a number }$\theta =\theta \left( \eta \right)
\in (0,1) $\emph{\ such that the following convergence estimate is valid:}%
\begin{equation*}
\left\Vert \widetilde{c}_{k,n}-c^{\ast }\right\Vert _{H^{1}\left(
\varepsilon ,b\right) }\leq C\theta ^{k}\left\Vert \widetilde{q}_{0,n}-%
\widetilde{q}_{n,\min }\right\Vert _{H^{4}\left( \Omega \right) }
\end{equation*}%
\begin{equation*}
+\frac{C^{n/2}}{\lambda ^{n/2}}\left( \dint\limits_{\Omega }\left( \overline{%
q}_{0x}^{2}+\overline{q}_{0t}^{2}+\overline{q}_{0}^{2}\right) \left(
x,t\right) e^{2\lambda \varphi }dxdt\right) ^{1/2}+C\left( \delta +\sqrt{%
\beta }\right) ,
\end{equation*}%
\emph{where functions }$\widetilde{c}_{k,n},\widetilde{q}_{0,n},\widetilde{q}%
_{n,\min }$\emph{\ and }$\overline{q}_{0}$\emph{\ are defined in (\ref{7}), (%
\ref{9}) and (\ref{10}) respectively. }

\section{Numerical Studies}

\label{sec10}

\subsection{Numerical implementation}

\label{sec10.1}

To generate the simulated data, by Lemma 2.2, we solve problem (\ref{2.3}), (%
\ref{2.4}) for the case when the whole real line is replaced by a large
interval $(-a,a)$ with the absorbing boundary conditions (\ref{2.21})--(\ref%
{2.22}). More precisely, just as in Section 6.1 of \cite{KLNSN}, we use the
implicit scheme to numerically solve 
\begin{equation}
\left\{ 
\begin{array}{rcll}
c(x)u_{tt}(x,t) & = & u_{xx}(x,t) & (x,t)\in (-a,a)\times (0,T), \\ 
u(-a,t)-u_{x}(-a,t) & = & 0 & t\in (0,T), \\ 
u(a,t)+u_{x}(a,t) & = & 0 & t\in (0,T), \\ 
u(x,0) & = & 0 & x\in \mathbb{R}, \\ 
u_{t}(x,0) & = & \widetilde{\delta }(x) & x\in \mathbb{R},%
\end{array}%
\right.  \label{6.1}
\end{equation}%
where $a=5$, $T=6$ and 
\begin{equation*}
\widetilde{\delta }(x)=\frac{30}{\sqrt{2\pi }}e^{-\frac{(30x)^{2}}{2}}
\end{equation*}%
is a smooth approximation of the Dirac function $\delta \left( x\right) $.
We solve problem (\ref{6.1}) by the implicit finite difference method. In
the finite difference scheme, we arrange a uniform partition for the
interval $[-a,a]$ as $\{y_{0}=-a,y_{1},\dots ,y_{N}=a\}\subset \lbrack -a,a]$
with $y_{i}=a+2ia/N_{x}$, $i=0,\dots ,N_{x}$, where $N_{x}$ is a large
number. In the time domain, we split the interval $[0,T]$ into $N_{t}+1$
uniform sub-intervals $[t_{j},t_{j+1}]$, $j=0,\dots ,N_{t}$, with $%
t_{j}=jT/N_{t},$ where $N_{t}$ is a large number. In our computational
setting, $N_{x}=3001$ and $N_{t}=301$. These numbers are the same as in \cite%
{KLNSN}.

We observe a computational error for the function $u$ near $(x=0,t=0)$. This
is due to the fact that the function $\widetilde{\delta }\left( x\right) $
is not exactly equal to the Dirac function. We correct the error as follows.
It follows from (\ref{2.2}), (\ref{2.6}) and (\ref{2.7}) that $u(x,t)=1/2$
in a neighborhood of the point $\left( x,t\right) =\left( 0,0\right) .$ We,
therefore, replace the data $u(\epsilon ,t)$ by $1/2$ when $|t|$ is small.
In our computation, we set $u(x,t)=1/2$ for $(x,t)\in \lbrack
0,0.0067]\times \lbrack 0,0.26]$. This data correction step is exactly the
same as in Section 6.1 of \cite{KLNSN} and is illustrated by Figure 2 in
that publication. Then, we can extract the noiseless data $g_{0}^{\ast }$
easily. We next add the noise into the data via the formula 
\begin{equation}
g_{0}=g_{0}^{\ast }(1+\delta \mbox{rand})  \label{10.2}
\end{equation}%
where $\delta $ is the noise level and rand is the function that generates
uniformly distributed random numbers in the range $[-1,1].$ In all numerical
tests with simulated data below, the noise level $\delta =0.05,$ i.e. $5\%.$
Due to \eqref{2.80}, the function $g_{1}=g_{0}^{\prime }.$ 
Due to the presence of noise, see \eqref{10.2}, we cannot compute $%
g_{1}=g_{0}^{\prime }$ by the finite difference method. Hence, the function $%
g_{0}^{\prime }$ is computed by the Tikhonov regularization method. The
version of the Tikhonov regularization method for this problem is
well-known. Hence, we do not describe this step here. 

Having the data for the function $q$ in hand, we proceed as in Algorithm \ref%
{alg}. 
\begin{algorithm}
	\caption{\label{alg} A numerical method to solve Problem 3.1}
	\begin{algorithmic}[1]
	\STATE \label{choice} Choose a set of  parameters $\lambda$, $\alpha$ and $\beta$. 
	\STATE  \label{choose q0} Compute the function $q_0$ by minimizing the functional $J^{(0)}_{\lambda, \beta}$ defined in \eqref{40.3}. 
	Due to \eqref{3.4}, the initial reconstruction is given by
	\begin{equation*}
		c_{\rm init}(x) = \frac{1}{(2 q_0(x, 0))^4} \quad \mbox{for all } x \in [\epsilon, b].
	\end{equation*}
	\STATE \label{step min} Assume that the function $q_{n - 1}$ is known. 
	We compute the function $q_n$ by minimizing the function $J^{(n)}_{\lambda, \beta}$ defined in \eqref{405}.
	\STATE  Set $q_{\rm comp} = q_n$ when $n = n^*$ is large enough.
	\STATE \label{step 4} Due to \eqref{3.4}, the function $c_{\rm comp}$ is set to be
	\[
		c_{\rm comp}(x) = \frac{1}{(2 q_{\rm comp}(x, 0))^4} \quad \mbox{for all } x \in [\epsilon, b].
	\]
 	\end{algorithmic}
\end{algorithm}In step \ref{choice} of Algorithm \ref{alg}, we choose $%
\lambda =2$, $\alpha =0.3$ and $\beta =10^{-11}.$ These parameters were
chosen by a trial-error process that is similar to the one in \cite{SKN1}.
Just as in \cite{SKN1}, we choose a reference numerical test in which we
know the true solution. In fact, test 1 of section 10.2 was our reference
test. We tried several values of $\lambda ,$ $\alpha ,$ and $\beta $ until
the numerical solution to that reference test became acceptable. Then, we
have used the same values of these parameters for all other tests, including
all five (5) available cases of experimental data.

We next implement Steps \ref{choose q0} and \ref{step min} of Algorithm \ref%
{alg}. We write differential operators in the functionals $J_{\lambda ,\beta
}^{(0)}$ and $J_{\lambda ,\beta }^{(n)}$ in the finite differences with the
step size in space $\Delta x = 0.0033$ and the step size in time $\Delta t =
0.02$ and minimize resulting functionals with respect to values of
corresponding functions at grid points. Since the integrand in the
definitions of the functional $J_{\lambda ,\beta }^{(n)},n=0,1,...$ is the
square of linear functions, then its minimizer is its critical point. 
In finite differences, we can write a linear system whose solution is the
desired critical point. We solve this system by the command
\textquotedblleft lsqlin" of Matlab. The details in implementation by the
finite difference method including the discretization, the derivation of the
linear system for the critical point, and the use of \textquotedblleft
lsqlin" are very similar to the scheme in \cite%
{Nguyen:CAMWA2020,Nguyens:jiip2020}. 
Recall that in our theory, in the definition of the functional $J_{\lambda
,\beta }^{(n)}$ acting on $q_{n}$, see \eqref{405}, we replaced $q_{n-1}$
with its analog $q_{n-1,\mathrm{min}}$ which belongs to the bounded set $%
\overline{B(R,g_{0})},$ and also replaced $q_{n-1,\mathrm{min}}(x,0)$ with $%
q_{n-1,\mathrm{min}}^{0}(x,0)$. These replacements are only for the
theoretical purpose to avoid the case when $q_{n-1}$ blows up. However, in
the numerical studies, these steps can be relaxed. This means that in Step %
\ref{step min}, we have minimized the finite difference analog of the
functional 
\begin{equation}
q_{n}\mapsto \int_{\Omega }\Big(q_{nxx}(x,t)-\frac{q_{nxt}(x,t)}{2\big(%
q_{n-1}(x,0)\big)^{2}}+\frac{\partial _{t}q_{n-1}\left( x,t\right) \partial
_{x}q_{n-1}(x,0)}{2\big(q_{n-1}(x,0)\big)^{3}}\Big)^{2}e^{2\lambda \varphi
}dxdt+\beta \left\Vert q_{n}\right\Vert _{H^{2}\left( \Omega \right) }^{2},
\label{10.3}
\end{equation}%
%
%
%
%
subject to the boundary conditions in lines 2 and 3 of (\ref{3.11}). Another
numerical simplification is that rather than using the $H^{4}-$norm in the
regularization term, we use the $H^{2}-$ norm in \eqref{10.3}. Although the
theoretical analysis supporting the above simplifications is missing, we did
not experience any difficulties in computing the numerical solutions to
Problem 3.1. All of our numerical results are satisfactory.


\subsection{Numerical results for computationally simulated data}

\label{sec10.2}

To test Algorithm \ref{alg}, we present four (4) numerical examples.

\textbf{Test 1 (the reference test).} We first test the case of one
inclusion with a high inclusion/background contrast. The true dielectric
constant function $c(x)$ has the following form 
\begin{equation}
c_{\mathrm{true}}(x)=\left\{ 
\begin{array}{ll}
1+14e^{\frac{(x-1)^{2}}{(x-1)^{2}-0.2^{2}}} & \mbox{if }|x-1|<0.2, \\ 
1 & \mbox{otherwise}.%
\end{array}%
\right.  \label{c1}
\end{equation}%
\begin{figure}[h]
\subfloat[\label{1a}]{\includegraphics[width = .45\textwidth]{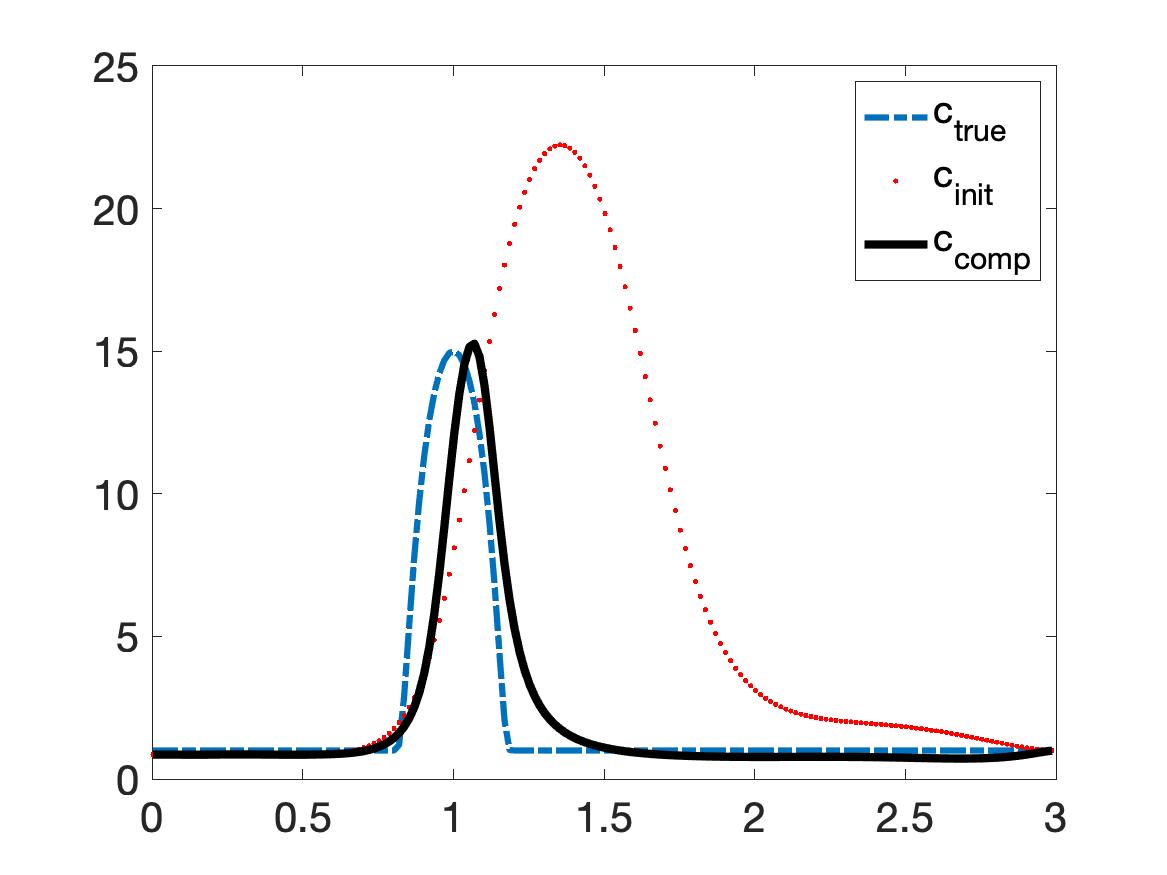}}
\quad \subfloat[\label{1b}]{\includegraphics[width = .45\textwidth]{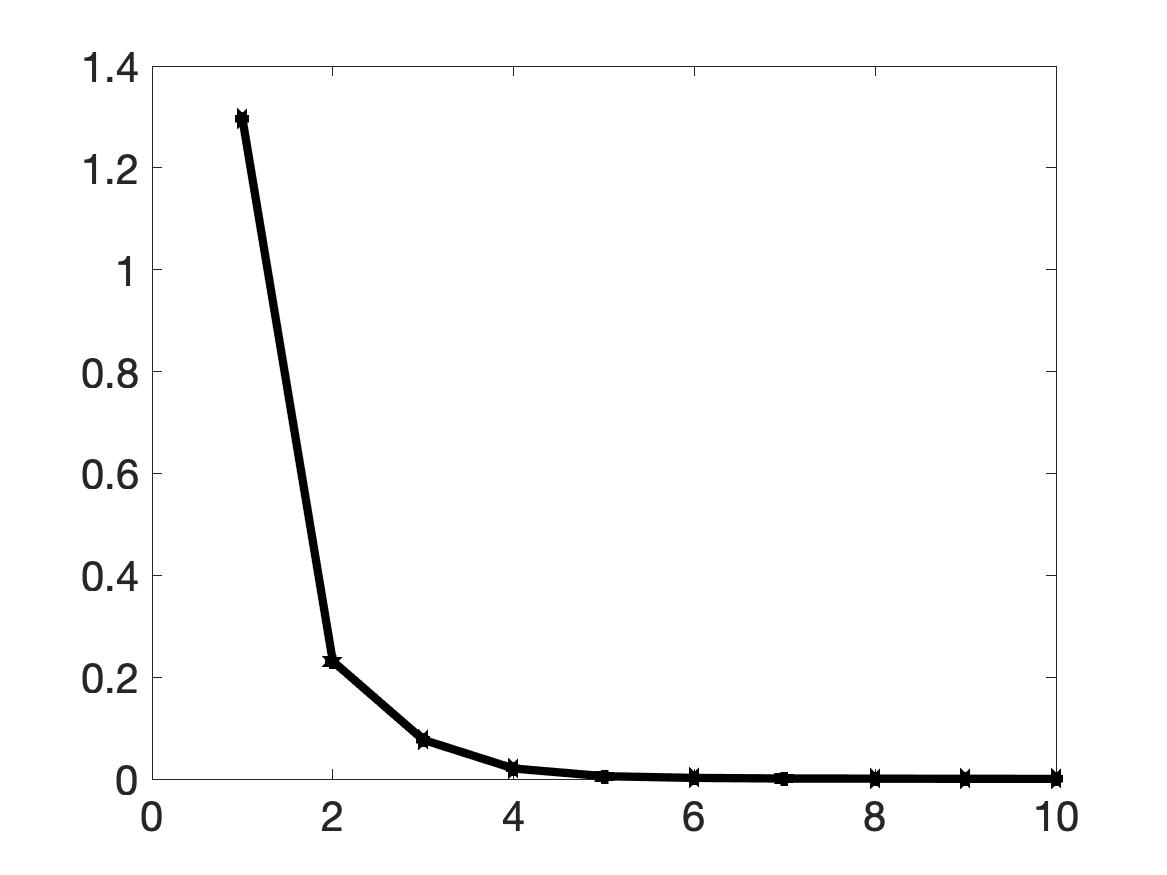}}
\caption{\textit{Test 1. The true and reconstructed functions $c(x),$ where $%
c_{\mathrm{true}}$ is given in \eqref{c1}. (a) The functions $c_{\mathrm{init%
}}$ and $c_{\mathrm{comp}}$ are obtained by Step \protect\ref{choose q0} and
Step \protect\ref{step 4} of Algorithm \protect\ref{alg} respectively. (b)
The consecutive relative error $\Vert c_{n}-c_{n-1}\Vert _{L^{\infty }(%
\protect\epsilon ,M)}/\Vert c_{n}\Vert _{L^{\infty }(\protect\epsilon ,M)}$, 
$n=1,\dots ,10.$ The data is with $\protect\delta =5\%$ noise. }}
\label{test1}
\end{figure}

Thus, the inclusion/background contrast as 15:1. It is evident from Figure %
\ref{test1} that we can successfully detect an object. The diameter of this
object is 0.4 and the distance between this object and the source is 1. The
true inclusion/background contrast is 15/1. The maximal value of the
computed dielectric constant is 15.28. The relative error in this maximal
value is 1.89\% while the noise level in the data is 5\%. Although the
contrast is high, our method provides good numerical solution without
requiring any knowledge of the true solution. Our method converges fast.
Although the initial reconstruction $c_{\mathrm{init}}$ computed in Step \ref%
{choose q0} of Algorithm \ref{alg} is not very good, see Figure \ref{1a},
one can see in Figure \ref{1b}, that the convergence occurs after only five
(5) iterations. This fact verifies numerically the \textquotedblleft
contraction property" of Theorem 8.1 including the key estimate \eqref{8.18}.

\begin{figure}[h!]
\subfloat[\label{2a}]{\includegraphics[width = .45\textwidth]{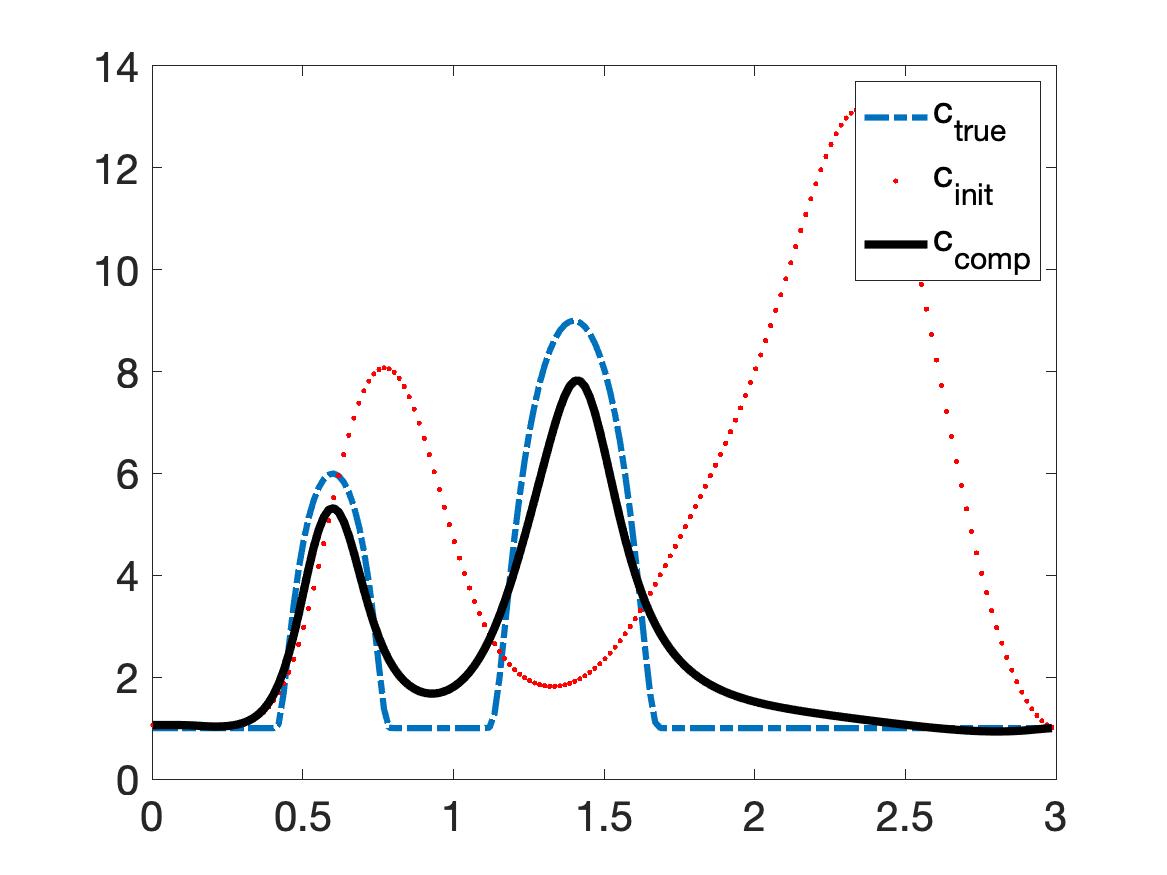}}
\quad \subfloat[\label{2b}]{\includegraphics[width = .45\textwidth]{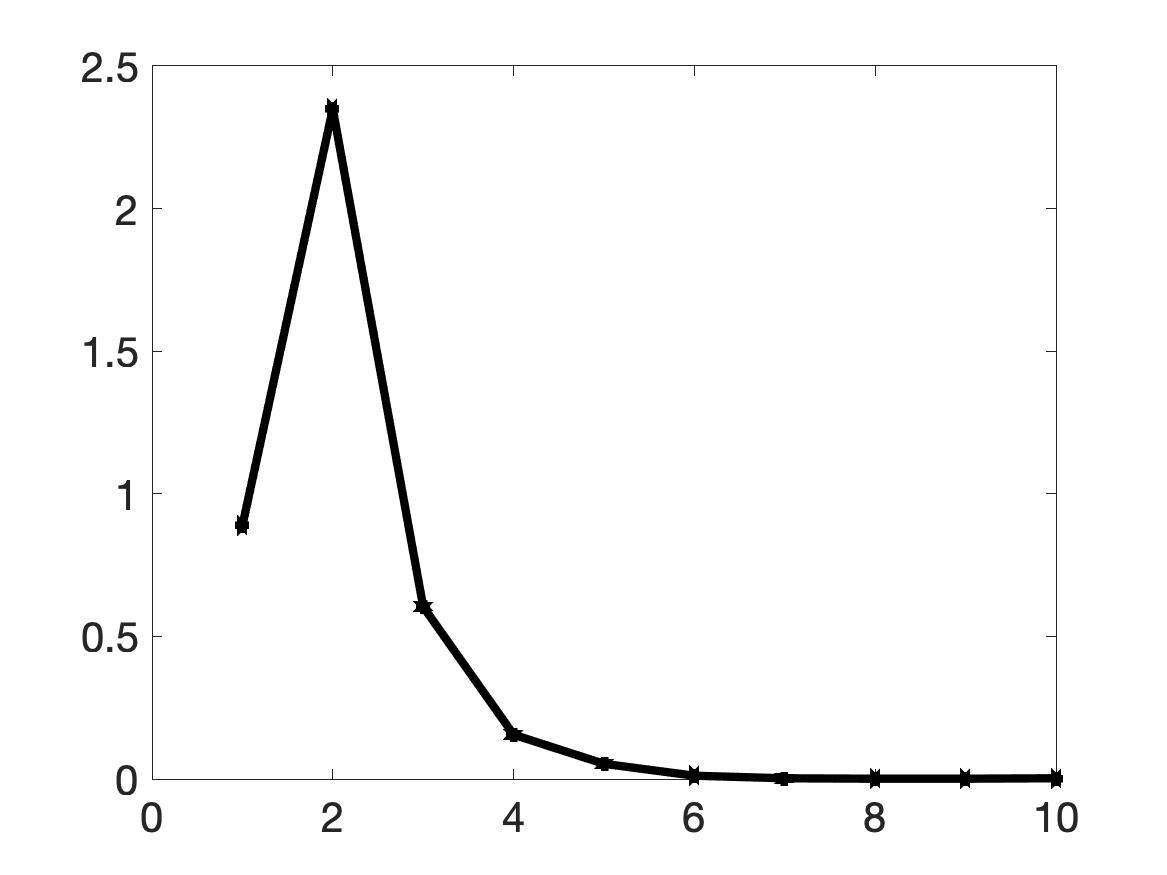}}
\caption{\textit{Test 2. The true and reconstructed functions $c(x),$ where $%
c_{\mathrm{true}}$ is given in \eqref{c2}. (a) The functions $c_{\mathrm{init%
}}$ and $c_{\mathrm{comp}}$ are obtained by Step \protect\ref{choose q0} and
Step \protect\ref{step 4} of Algorithm \protect\ref{alg} respectively. (b)
The consecutive relative error $\|c_{n} - c_{n-1}\|_{L^{\infty}(\protect%
\epsilon, M)}/\|c_{n} \|_{L^{\infty}(\protect\epsilon, M)}$, $n = 1, \dots,
10.$ The data is with $\protect\delta = 5\%$ noise. }}
\label{test2}
\end{figure}

\textbf{Test 2.} The true function $c$ in this test has two
\textquotedblleft inclusions", 
\begin{equation}
c_{\mathrm{true}}(x)=\left\{ 
\begin{array}{ll}
1+5e^{\frac{(x-0.6)^{2}}{(x-0.6)^{2}-0.2^{2}}} & \mbox{ if }\left\vert
x-0.6\right\vert <0.2, \\ 
1+8e^{\frac{(x-1.4)^{2}}{(x-1.4)^{2}-0.3^{2}}} & \mbox{ if }\left\vert
x-1.4\right\vert <0.3, \\ 
1 & \mbox{ otherwise. }%
\end{array}%
\right.  \label{c2}
\end{equation}%
Numerical results of this test are displayed in Figure \ref{test2}.

Test 2 is more complicated than Test 1. However, we still obtain good
numerical results. It is evident from Figure \ref{2a} that we can precisely
detect the two inclusions without any initial guess. The true maximal values
of the dielectric constant of the inclusions in the left and the right are 6
and 9 respectively. The reconstructed maximal values in those inclusions are
acceptable. They are 5.31 (relative error 11.5\%) and 7.8 (relative error
13.3\%). As in Test 1, the initial reconstruction $c_{\mathrm{init}}$
computed in Step \ref{choose q0} of Algorithm \ref{alg} is far from $c_{%
\mathrm{true}}(x).$ Still, our iterative procedure converges after 7
iterations, see Figure \ref{2b}.

\begin{figure}[h]
\subfloat[\label{3a}]{\includegraphics[width = .45\textwidth]{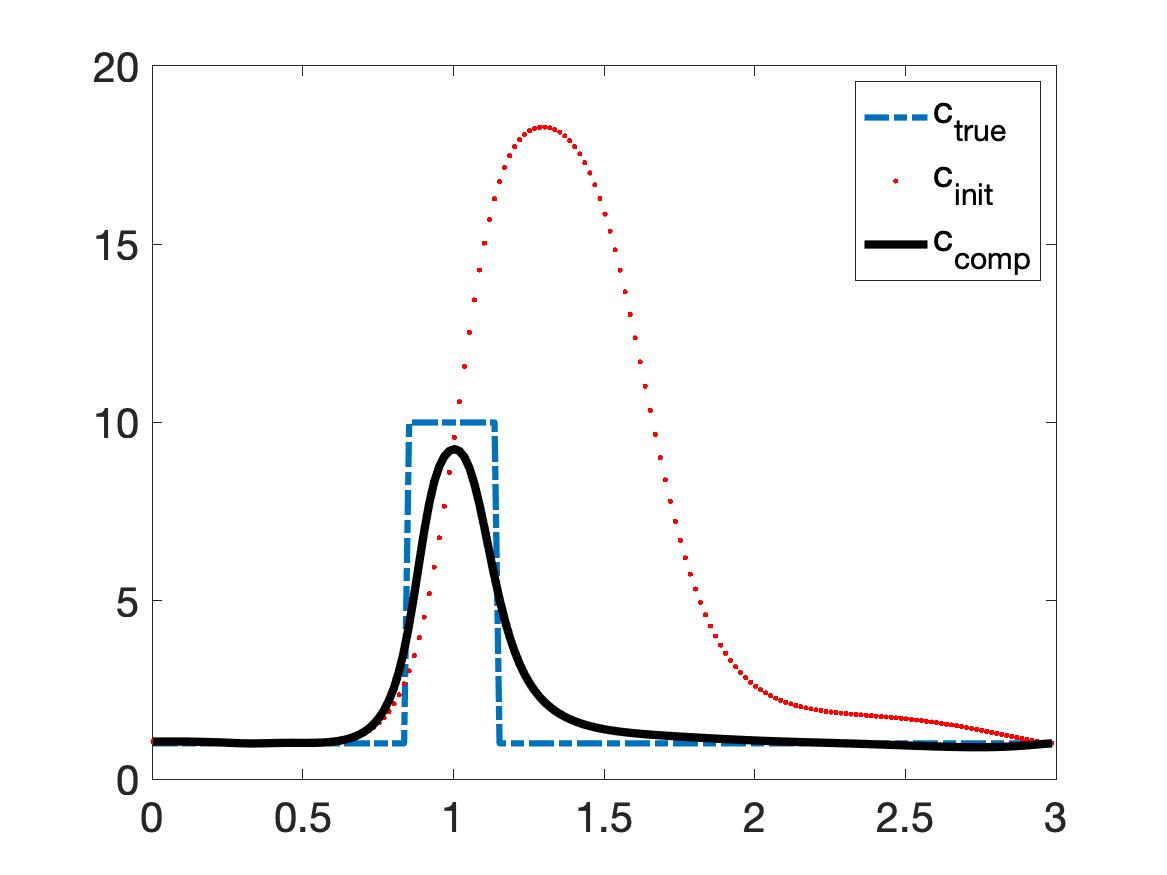}}
\quad \subfloat[\label{3b}]{\includegraphics[width = .45\textwidth]{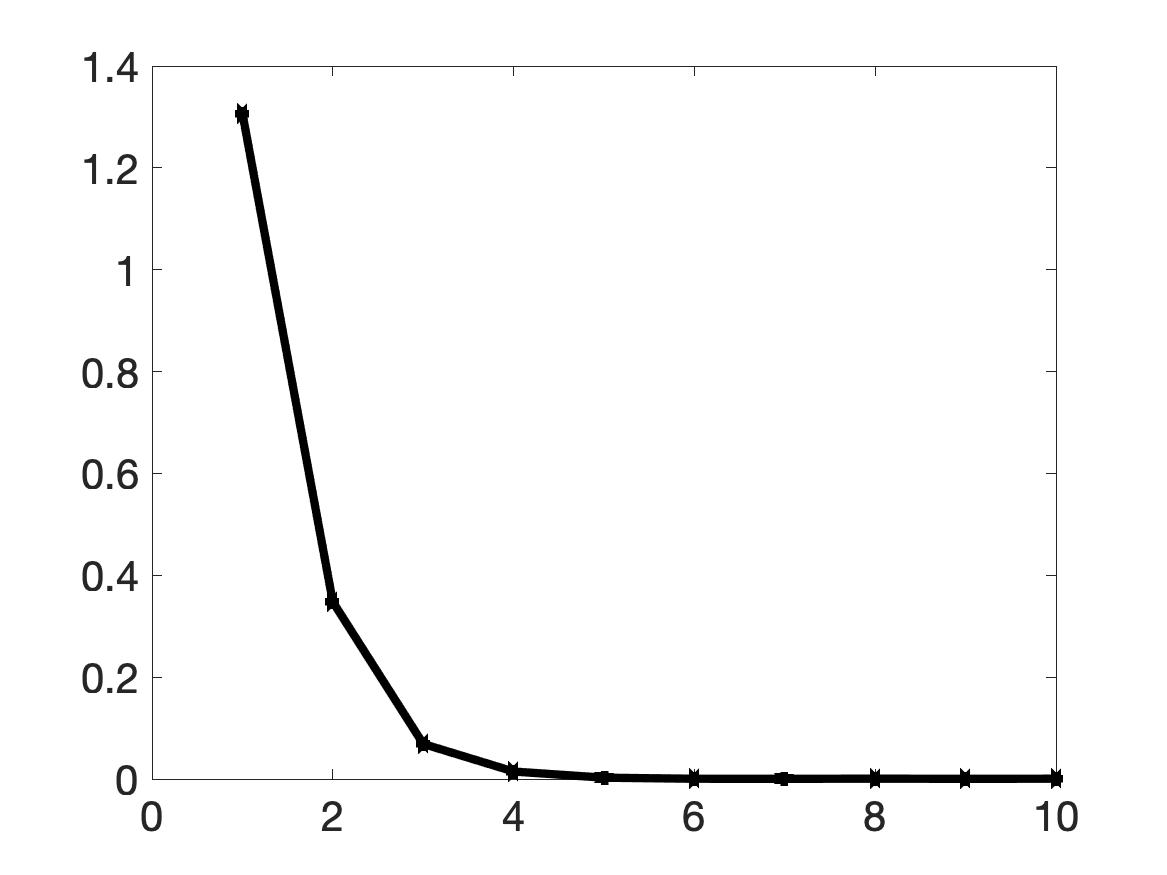}}
\caption{\textit{Test 3. The true and reconstructed functions $c(x),$ where $%
c_{\mathrm{true}}$ is given in \eqref{c3}. (a) The functions $c_{\mathrm{init%
}}$ and $c_{\mathrm{comp}}$ are obtained by Step \protect\ref{choose q0} and
Step \protect\ref{step 4} of Algorithm \protect\ref{alg} respectively. (b)
The consecutive relative error $\Vert c_{n}-c_{n-1}\Vert _{L^{\infty }(%
\protect\epsilon ,M)}/\Vert c_{n}\Vert _{L^{\infty }(\protect\epsilon ,M)}$, 
$n=1,\dots ,10.$ The data is with $\protect\delta =5\%$ noise. }}
\label{test3}
\end{figure}

\textbf{Test 3} We test the case when the function $c_{\mathrm{true}}(x)$ is
discontinuous. It is a piecewise constant function given by 
\begin{equation}
c_{\mathrm{true}}(x)=\left\{ 
\begin{array}{ll}
10 & \mbox{ if }\left\vert x-1\right\vert <0.15, \\ 
1 & \mbox { otherwise. }%
\end{array}%
\right.  \label{c3}
\end{equation}%
The numerical solution for this test is presented in Figure \ref{test3}. 

Although the function $c_{\text{true}}$ is not smooth and actually not even
continuous, Algorithm \ref{alg} works and provides a reliable numerical
solution. The computed maximal value of the dielectric constant of the
object is 9.25 (relative error 7.5\%), which is acceptable. The location of
the object is detected precisely, see Figure \ref{3a}. As in the previous
two tests, Algorithm \ref{alg} converges fast. After the fifth iteration,
the reconstructed function $c_{n}$ becomes unchanged. Again, this fact
numerically confirms both Theorem 8.1 and the robustness of our method.

\begin{figure}[h]
\subfloat[\label{5a}]{\includegraphics[width = .45\textwidth]{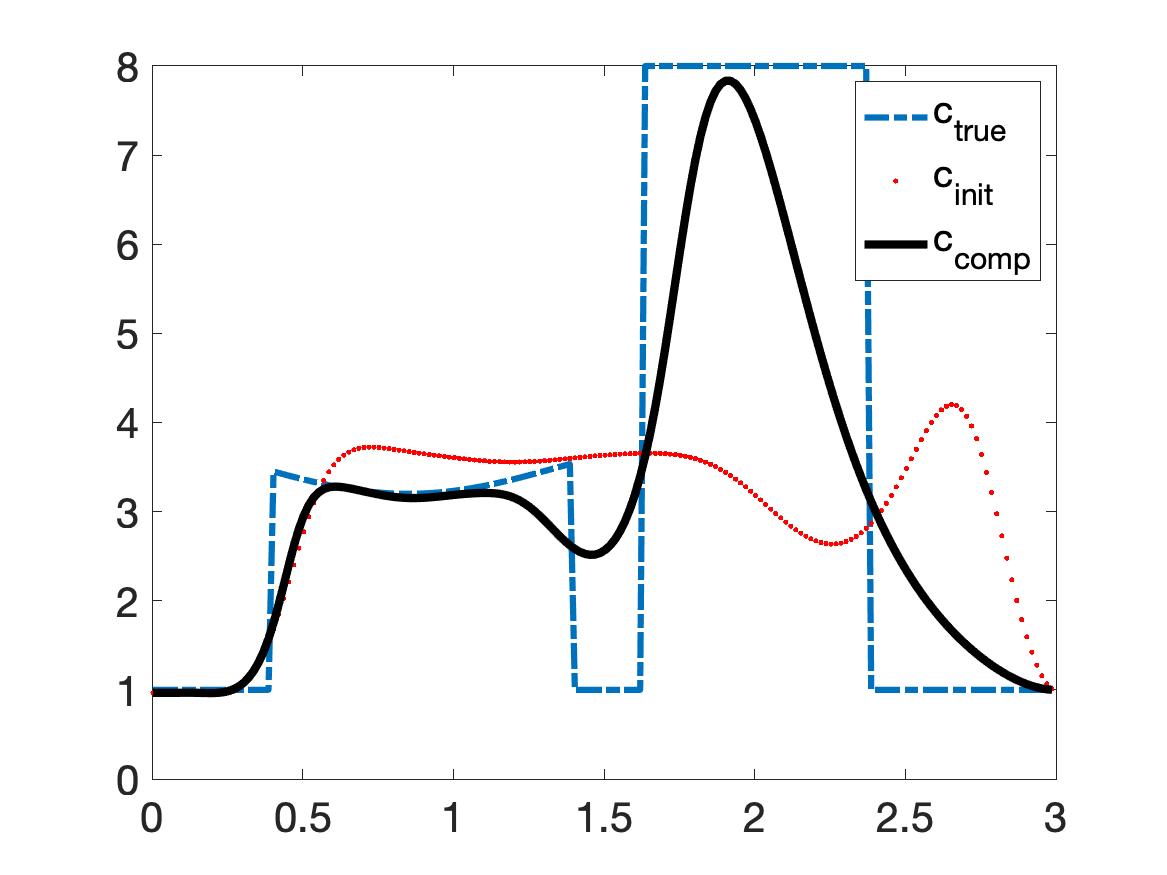}}
\quad \subfloat[\label{5b}]{\includegraphics[width = .45\textwidth]{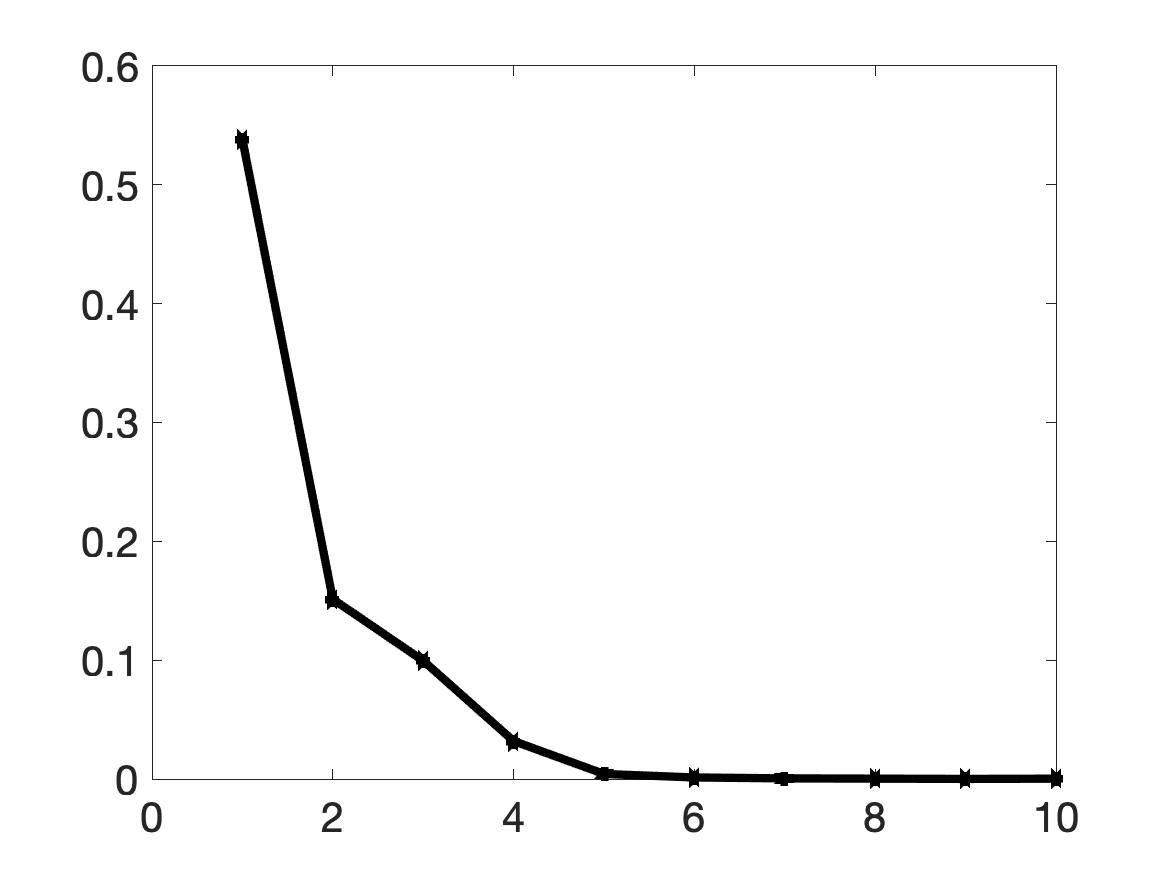}}
\caption{\textit{Test 4. The true and reconstructed functions $c(x),$ where $%
c_{\mathrm{true}}$ is given in \eqref{c5}. (a) The functions $c_{\mathrm{init%
}}$ and $c_{\mathrm{comp}}$ are obtained by Step \protect\ref{choose q0} and
Step \protect\ref{step 4} of Algorithm \protect\ref{alg} respectively. (b)
The consecutive relative error $\Vert c_{n}-c_{n-1}\Vert _{L^{\infty }(%
\protect\epsilon ,M)}/\Vert c_{n}\Vert _{L^{\infty }(\protect\epsilon ,M)}$, 
$n=1,\dots ,10.$ The data is with $\protect\delta =5\%$ noise. }}
\label{test5}
\end{figure}

\textbf{Test 4} In this test, the function $c_{\mathrm{true}}(x)$ has the
following form: 
\begin{equation}
c_{\mathrm{true}}(x)=\left\{ 
\begin{array}{ll}
3.5+0.3\ast \sin (\pi (x-1.35)) & \mbox{ if }\left\vert x-0.9\right\vert
<0.5, \\ 
8 & \mbox{ if }\left\vert x-2\right\vert <0.37, \\ 
1 & \mbox { otherwise. }%
\end{array}%
\right.  \label{c5}
\end{equation}%
This test is interesting since the graph of the function \eqref{c5} consists
of a curve and an inclusion. The numerical solution for this case is
presented in Figure \ref{test5}. 

One can observe from Figure \ref{5a} that our method to compute the initial
reconstruction in Step \ref{choose q0} of Algorithm \ref{alg} is not very
effective. However, after only 6 iterations, good numerical results are
obtained. The curve in the first inclusion locally coincides with the true
one and the maximal value of the computed dielectric constant within
inclusion is quite accurate: it is 7.83 (relative error 2.12\%). Our method
converges at the iteration number 6.

\begin{remark}
It follows from all above tests that Algorithm \ref{alg} is robust in
solving a highly nonlinear and severely ill-posed Problem 3.1. It provides
satisfactory numerical results with a fast convergence. For each test, the
computational time to compute the numerical solution is about 29 seconds on
a MacBook Pro 2019 with 2.6 GHz processor and 6 Intel i7 cores. This is
almost a real time computation.
\end{remark}

\section{Numerical Results for Experimental Data}

\label{sec11}

The experimental data we use to test Algorithm \ref{alg} were collected by
the Forward Looking Radar built in the US Army Research Laboratory \cite%
{NguyenWong}. The data were initially collected to detect and identify
targets mimicking shallow anti-personnel land mines and IEDs. Five mimicking
devices were: a bush, a wood stake, a metal box, a metal cylinder, and a
plastic cylinder. The bush and the wood stake were placed in the air, while
the other three objects were buried at a few centimeters depth in the
ground. We refer the reader to \cite{NguyenWong} for details about the
experimental set up and the data collection. This was recently repeated in 
\cite{KLNSN}. We, therefore, do not describe here the device and the data
collection procedure again. Since the location of the targets can be
detected by the Ground Position System (GPS), we are only interested in
computing the values of the dielectric constants of those targets. We do so
by our new method written in Algorithm \ref{alg}.

As in the earlier works using these data \cite{Karch, KLNSN, KlibLoc,
KlibKol1, KlibKol2, Kuzh1,Kuzh2, SKN1}, we compute the function $c_{\mathrm{%
rel}}(x),$ where 
\begin{equation}
c_{\mathrm{rel}}(x)=\left\{ 
\begin{array}{ll}
\frac{c_{\mathrm{target}}}{c_{\mathrm{bckgr}}}\left( x\right) \text{ if }%
\max \frac{c_{\mathrm{target}}}{c_{\mathrm{bckgr}}}\left( x\right) >1 & 
\text{and }x\in D, \\ 
1 & \text{otherwise},%
\end{array}%
\right.  \label{11.11}
\end{equation}%
\begin{equation}
c_{\mathrm{rel}}(x)=\left\{ 
\begin{array}{ll}
\min \frac{c_{\mathrm{target}}}{c_{\mathrm{bckgr}}}\left( x\right) \text{ if 
}\max \frac{c_{\mathrm{target}}}{c_{\mathrm{bckgr}}}\left( x\right) \leq 1 & 
\text{and }x\in D, \\ 
1 & \text{otherwise},%
\end{array}%
\right.  \label{11.12}
\end{equation}%
where $D$ is a sub interval of the interval $[\epsilon ,M],$ which is
occupied by the target. Next, we define the computed value of $c_{\mathrm{%
target}}$ as \cite{KlibLoc}: 
\begin{equation}
c_{\mathrm{comp}}=c_{\mathrm{bckgr}}\times \left\{ 
\begin{array}{c}
\max c_{\mathrm{rel}}(x)\text{ if }\max c_{\mathrm{rel}}(x)>1, \\ 
\min c_{\mathrm{rel}}(x)\text{ if }\max c_{\mathrm{rel}}(x)<1.%
\end{array}%
\right.  \label{11.1}
\end{equation}

As in the above cited publications, we have to preprocess the raw data of 
\cite{NguyenWong} before importing them into our solver. The data
preprocessing is exactly the same as in \cite[Section 7.1]{KLNSN}. First, we
observe that the $L_{\infty }-$norm of the experimental data far exceeds the 
$L_{\infty }-$norm of the computationally simulated data. Therefore, we need
to scale the experimental data by dividing it by a calibration factor $\mu
>0,$ i.e. we replace the raw experimental data $f_{\text{raw}}\left(
t\right) $ with  $f_{\text{scale}}\left( t\right) =f_{\text{raw}}\left(
t\right) /\mu $. Then we work only with  $f_{\text{scale}}\left( t\right) .$
To find the calibration factor, we use a trial-and-error process. First, we
select a reference target. We then try many values of $\mu $ such that the
reconstruction of the reference target is satisfactory, i.e. the computed
dielectric constant is in its published range, see below in this section
about the published range. Then this calibration factor is used is all other
tests. When the object is in the air, our reference target is bush. In this
case, the calibration factor $\mu _{\mathrm{air}}=534,592.$ When the object
is buried under the ground, our reference target is the metal box and the
calibration factor was $\mu _{\mathrm{ground}}=265,223$. 

Next, we preprocess the data $f_{\text{scale}}\left( t\right) $ as follows.
First, we first use a lower envelop (built in Matlab) to bound the data. We
then truncate the data that is not in a small interval centered at the
global minimizer of the data, see \cite[Section 7.1]{KLNSN} for the choice
of this small interval. But in the case of the plastic cylinder we use the
upper envelop. The choice of the upper or lower envelopes is as follows. We
look at the function $f_{\text{scale}}\left( t\right) $ and find the three
extrema with largest absolute values. If the middle extremal value among
these three is a minimum, then we bound the data by a lower envelop.
Otherwise, we use an upper envelop. See \cite[Section 7.1]{KLNSN} for more
details and the reason of this choice. In particular, Figures 4a, 4b, 4d,
5a, 5b, 5d, 5e, 5g, 5h of \cite{KLNSN} provide illustrations. Likewise, our
Figures 5 displays computed functions $c_{\text{target}}\left( x\right) $
for our five targets. The computed dielectric constants $c_{\mathrm{comp}}$
defined in \eqref{11.1} by Algorithm \ref{alg} is listed in Table \ref{tab1}.

\begin{figure}[h!]
\subfloat[$c_{\rm target}$ for
bush]{\includegraphics[width=.3\textwidth]{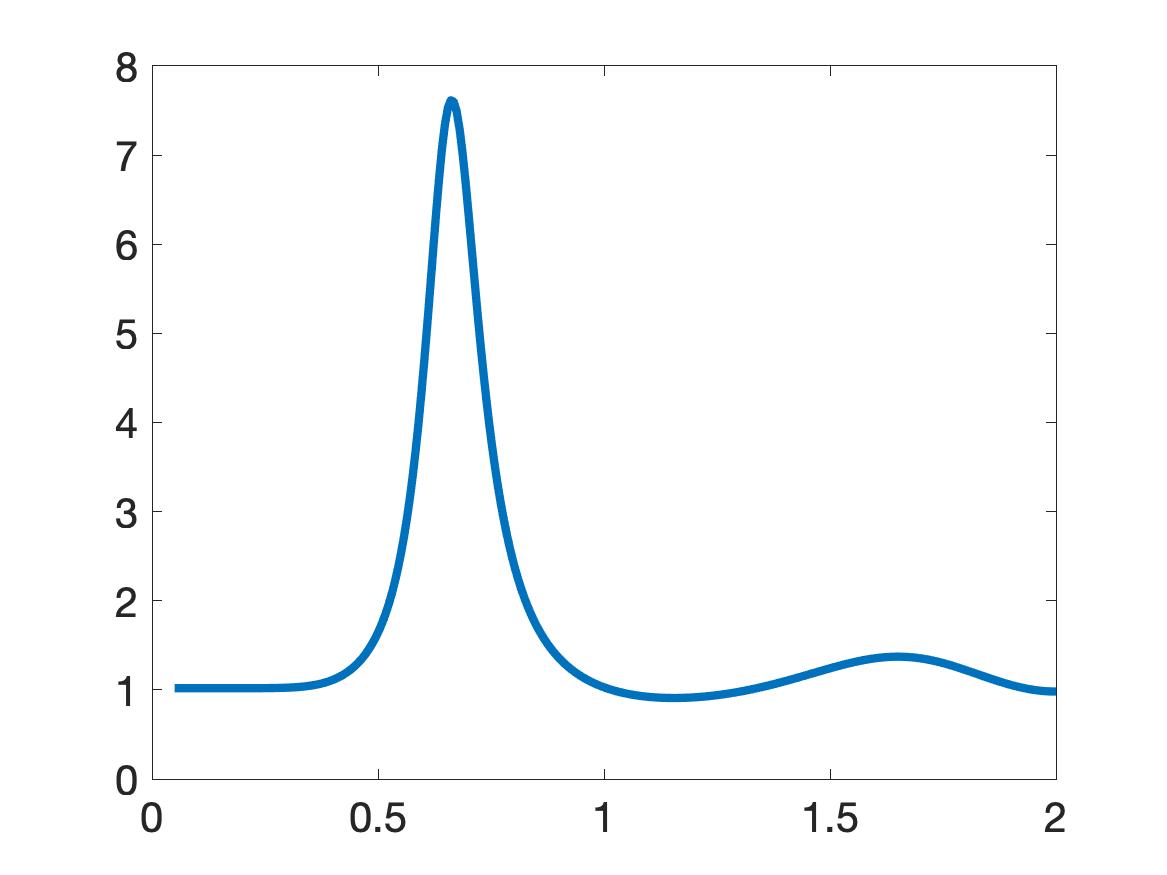}} \quad 
\subfloat[$c_{\rm target}$ for wood
stake]{\includegraphics[width=.3\textwidth]{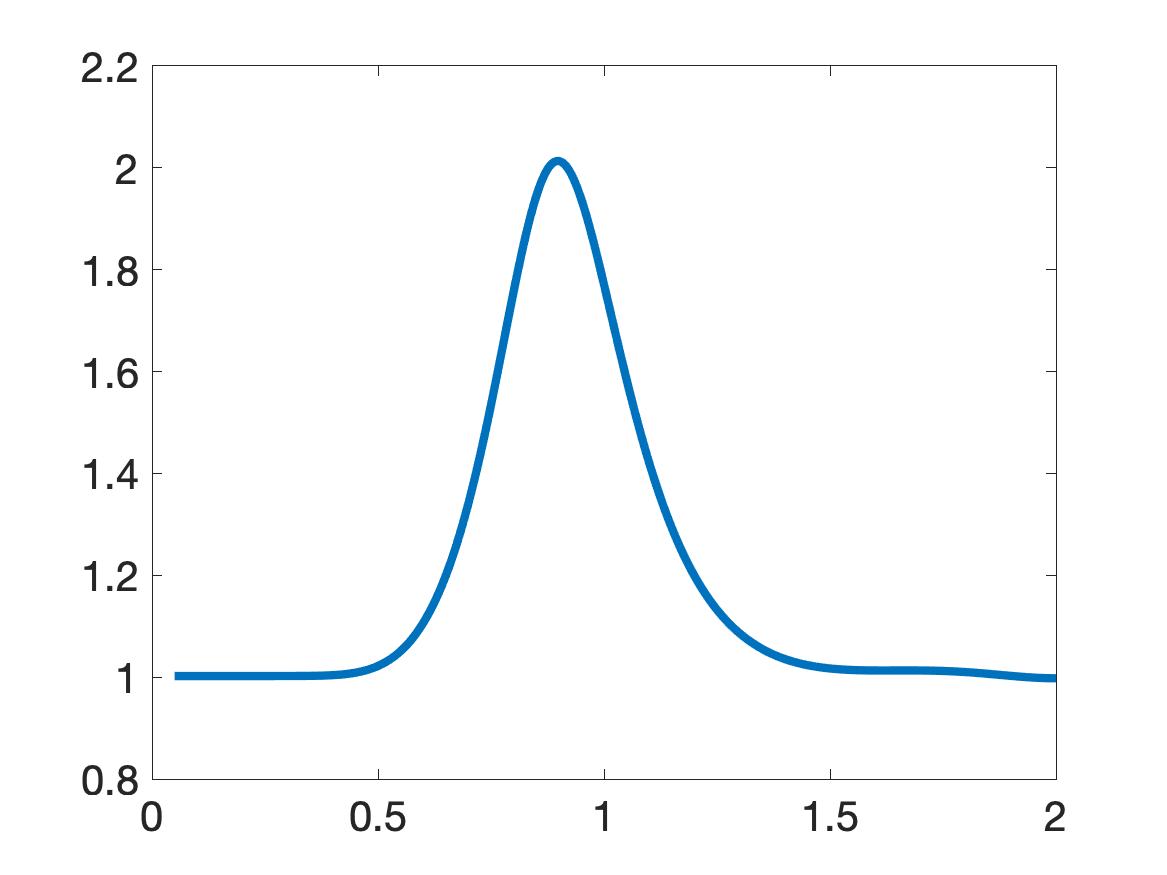}}
\par
\subfloat[$c_{\rm target}$ for metal
box]{\includegraphics[width=.3\textwidth]{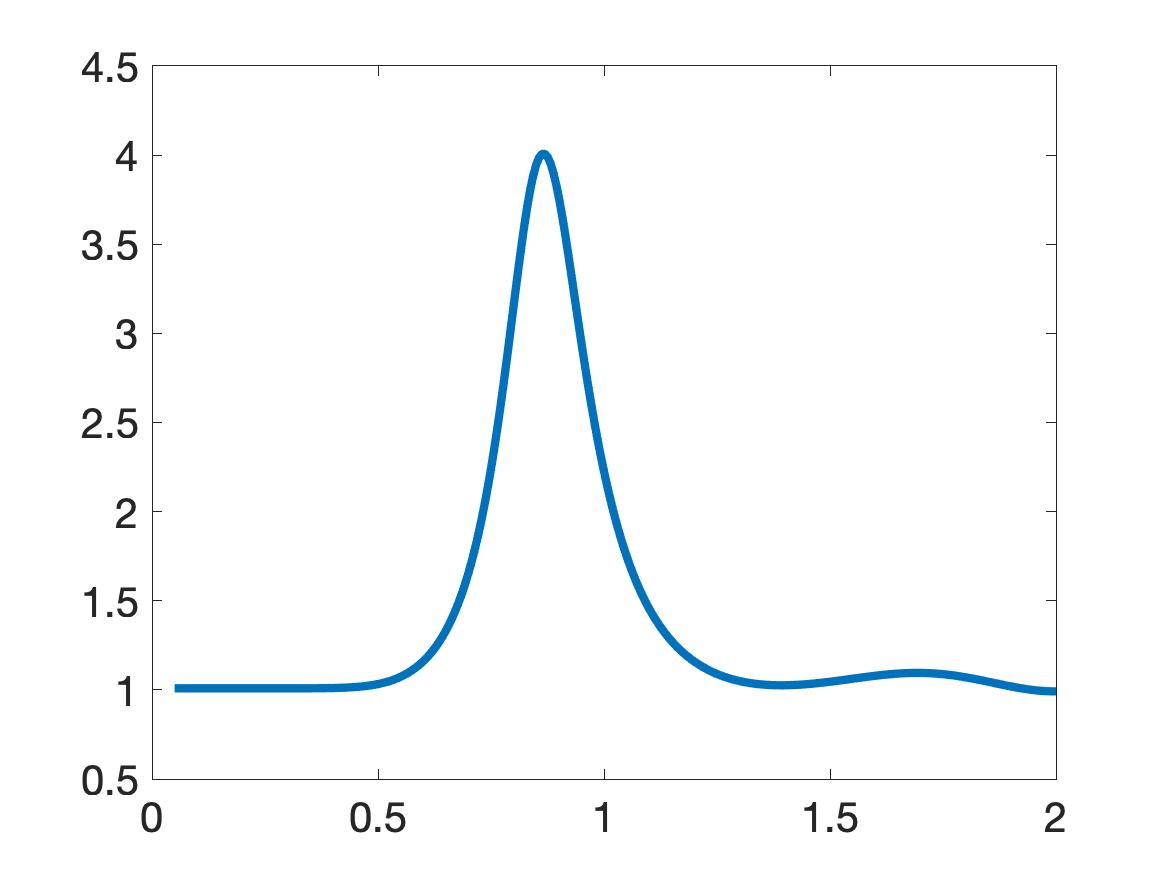}} \quad 
\subfloat[$c_{\rm target}$ for metal
cylinder]{\includegraphics[width=.3\textwidth]{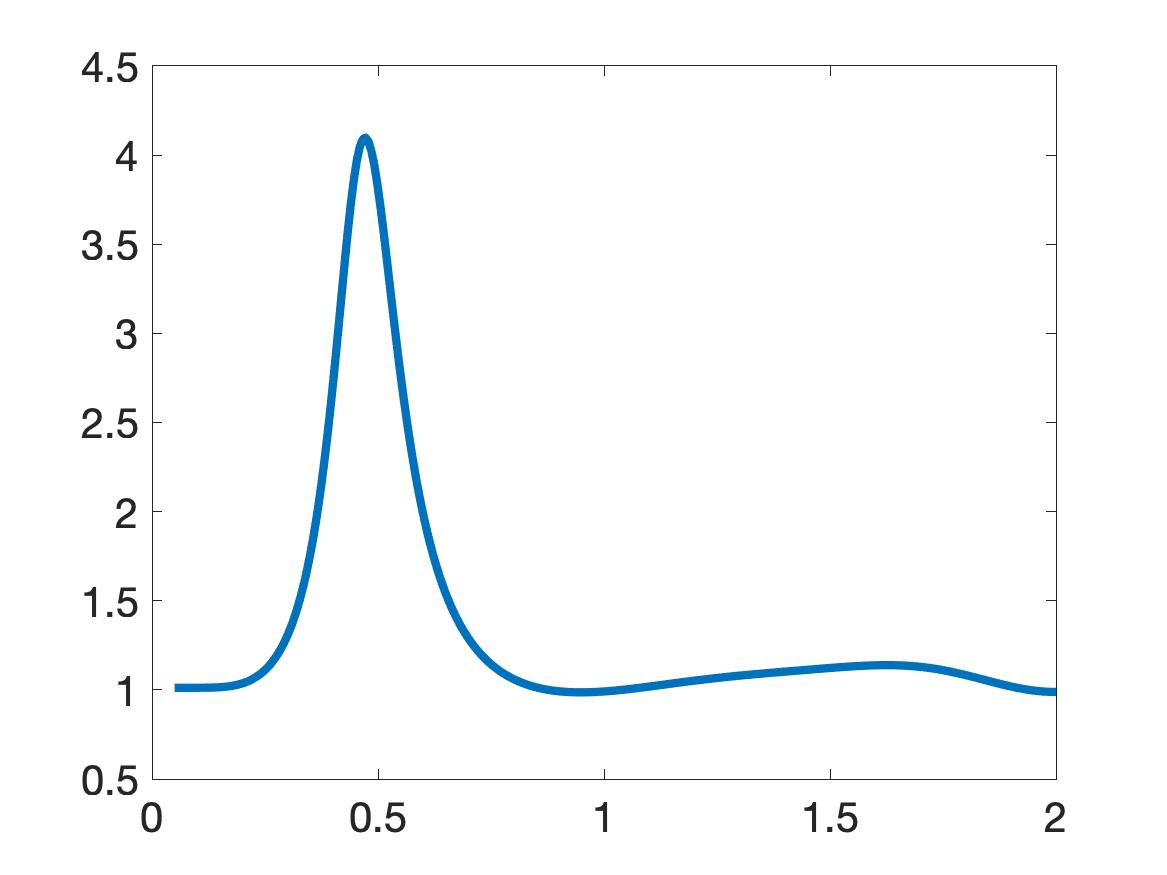}} \quad 
\subfloat[$c_{\rm target}$ for plastic
cylinder]{\includegraphics[width=.3\textwidth]{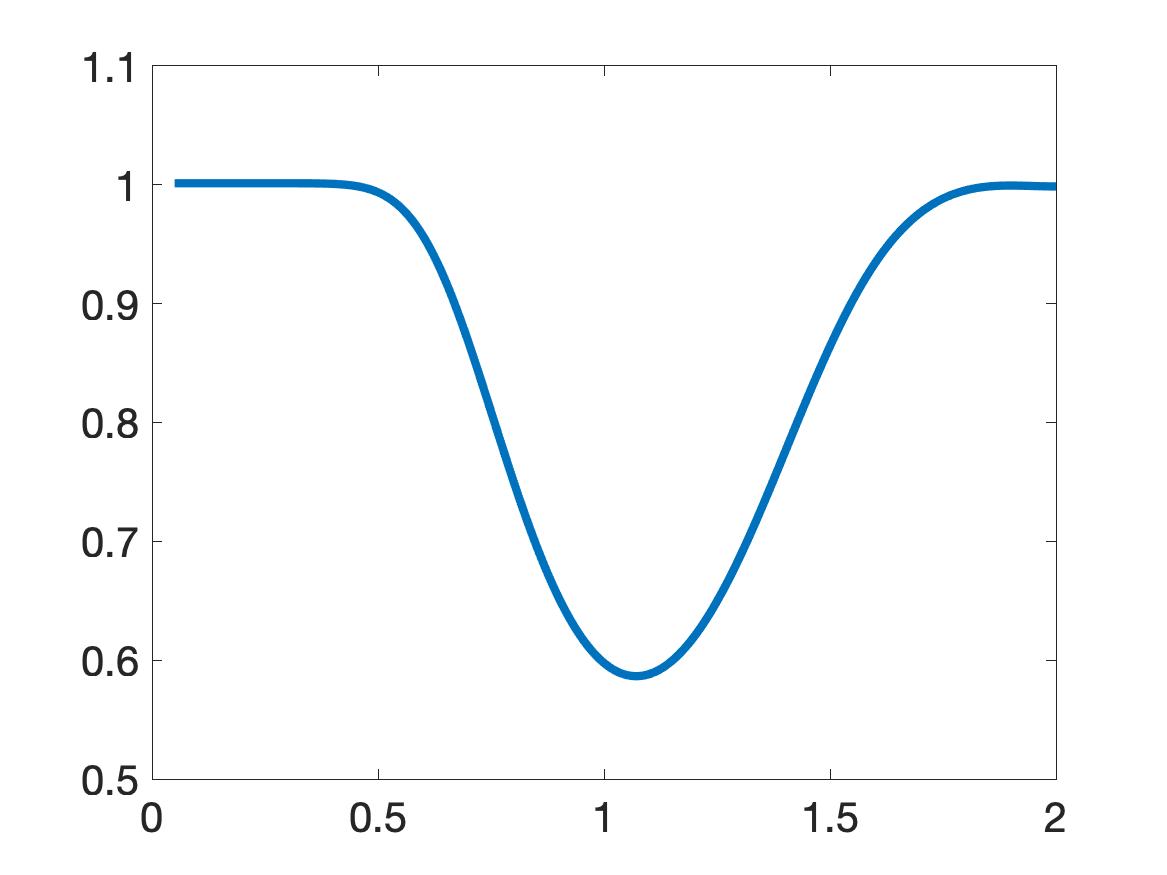}}
\caption{\textit{Computed functions $c_{\mathrm{target}}(x, y)$ for our five
targets, also see (\protect\ref{11.11})- (\protect\ref{11.1}) and Table 
\protect\ref{tab1}.}}
\end{figure}

\begin{table}[h!]
\begin{center}
\begin{tabular}{|c|c|c|c|c|c|}
\hline
Target & $c_{\mathrm{bckgr}}$ & computed $c_{\mathrm{rel}}$ & $c_{\mathrm{%
bckgr}}$ & computed $c_{\text{target}}$ & True $c_{\text{target}}$ \\ \hline
Bush & 1 & 7.62 & 1 & 10.99 & $[3,20]$ \\ 
Wood stake & 1 & 2.01 & 1 & 2.26 & $[2,6]$ \\ 
Metal box & 4 & 4.00 & $[3,5]$ & $[12.00,20.00]$ & $[10,30]$ \\ 
Metal cylinder & 4 & 4.01 & $[3,5]$ & $[12.3,20.5]$ & $[10,30]$ \\ 
Plastic cylinder & 4 & 0.59 & $[3,5]$ & $[1.6,2.95]$ & $\left[ 1.1,3.2\right]
$ \\ \hline
\end{tabular}%
\\[0pt]
\end{center}
\caption{Computed dielectric constants of five targets}
\label{tab1}
\end{table}

The true values of dielectric constants of our targets were not measured in
the experiment. Therefore, we compare our computed values with the published
ones. The published values of the dielectric constants of our targets are
listed in the last column of Table \ref{tab1}. They can be found on the
website of Honeywell (Table of dielectric constants, https://goo.gl/kAxtzB).
Also, see \cite{Chuah} for the experimentally measured range of the
dielectric constants of vegetation samples, which we assume have the same
range as the dielectric constant of bush. In the table of dielectric
constants of Honeywell as well as in \cite{Chuah}, any dielectric constant
is not a number. Rather, each dielectric constant of these references is
given within a certain interval. As to the metallic targets, it was
established in \cite{Kuzh1} that they have the so-called \textquotedblleft
apparent" dielectric constant whose values are in the interval $\left[ 10,30%
\right] .$

\textbf{Conclusion}. It is clear from Table \ref{tab1} that our computed
dielectric constants are consistent with the intervals of published ones.
Therefore, our results for experimental data are satisfactory.

\begin{remark}[Speed of computations]
Our experimental data are sparse. The size of the data in time is $N_{t}=80$%
, which is a lot smaller than that in the dense simulated data ($N_{t}=300$%
). Therefore, the speed of computations is much faster than for the case of
simulated data of section 10. Thus, all results of this section were
computed in real time on the same computer (MacBook Pro 2019 with 2.6 GHz
processor and 6 Intel i7 cores).
\end{remark}

\section{Concluding Remarks}

\label{sec12}

We have developed a new globally convergent numerical method for a 1-D
Coefficient Inverse Problem with backscattering data for the wave-like PDE (%
\ref{2.3}). This is the second generation of the above cited convexification
method of our research group. The main novelty here is that, rather than
minimizing a globally strictly convex weighted cost functional arising in
the convexification, we solve on each iterative step a linear boundary value
problem. This is done using the so-called Carleman Quasi-Reversibility
Method. Just like in the convexification, the key element of the technique
of this paper, on which our convergence analysis is based, is the presence
of the Carleman Weight Function in each quadratic functional which we
minimize. The convergence estimate is similar to the key estimate of the
classical contraction mapping principle. The latter explains the title of
this paper. We have proven a global convergence theorem of our method. Our
numerical results for computationally simulated data demonstrate high
reconstruction accuracies in the presence of 5\% random noise in the data.

Furthermore, our numerical results for experimentally collected data have
satisfactory accuracy via providing values of computed dielectric constants
of explosive-like targets within their published ranges.

A practically important observation here is that our computations for
experimental data were performed in \emph{real time}.

\section*{Acknowledgement}

The effort of TTL, MVK and LHN was supported by US Army Research Laboratory
and US Army Research Office grant W911NF-19-1-0044.


\end{document}